# Characteristic numbers and characteristic equations of parity vectors of Collatz sequences

Raouf Rajab
raouf.rajab@enig.rnu.tn

**ABSTRACT**

The present work deals with the characterization of parity vectors of Collatz sequences (of finite and infinite length). Such a characterization leads to the determination of several numbers (integers or non-integers) that we call the characteristic numbers of a given parity vector. Some characteristic numbers are linked together by equations that can be called characteristic equations of the considered parity vector. If a parity vector v of finite length n contains the first n terms of a parity vector V of infinite length then all the characteristic numbers of v are considered as characteristic numbers of order n of the infinite vector V. The limits of $n^{th}$ order characteristic numbers when n tends to infinity constitute the absolute characteristic numbers of the infinite vector V and they allow determining its behavior and properties. In this paper, we study the properties of parity vectors of infinite length based on their characteristic numbers. Then, we establish the relations between the first term of a Collatz sequence and its parity vector. Finally, still based on the characteristic numbers, we determine some conditions of existence and non-existence of divergent sequences.

**RESUME**

Le présent travail porte sur la caractérisation des vecteurs de parité des suites de Collatz (de longueurs finies et de longueurs infinies). Une telle caractérisation conduit à la détermination des plusieurs nombres (entiers ou non entiers) qu'on les appelle les nombres caractéristiques d'un vecteur de parité donné. Certain nombres caractéristiques sont liés entre eux par des équations qu'on peut les appeler équations caractéristiques du vecteur de parité considéré. Si un vecteur de parité v de longueur finie n renferme les n premiers termes d'un vecteur de parité V de longueur infinie alors tous les nombres caractéristiques de v sont considéré comme des nombres caractéristiques d'ordre n du vecteur infini V. Les limites de nombres caractéristiques d'ordre n lorsque n tend vers l'infini constituent les nombres caractéristiques absolus du vecteur infini V et ils permettent de déterminer ses comportements et ses propriétés. Dans cet article, on étudie les propriétés de vecteurs de parité de longueurs infinies en se basant sur leurs nombres caractéristiques. Ensuite, on établi les relations entre le premier terme d'une suite de Collatz et de son vecteur de parité. Enfin, en se basant toujours sur les nombres caractéristiques, on détermine quelques conditions d'existence et de non existence des suites divergentes.

**Mots clé :** Suites de Collatz ; Vecteur de parité ; Nombres caractéristiques ; Equations caractéristiques.

# 1. Introduction et préliminaire

En mécanique des fluides (ou en transfert), des centaines des nombres adimensionnels comme par exemple le nombre de Reynolds (Re), de Prandtl (Pr), de Nusselt,… sont utilisées pour étudier, caractériser et interpréter les phénomènes d'écoulement des fluides et les phénomènes des transferts au sein de ces fluides. Ces nombres contribuent d'une façon générale dans la compréhension des comportements des fluides en écoulements quelque soit leurs natures.

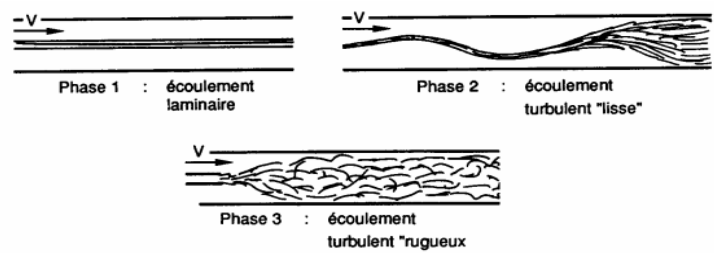

**FiG 1 :** Régimes d'écoulement et nombre de Reynolds

D'une manière similaire, dans le présent travail, on va explorer une nouvelle piste pour l'étude et l'analyse du problème 3n+1 basée sur la caractérisation des vecteurs de parité des suites de Collatz qui fait appelle à des nombres caractéristiques intervenant dans la détermination des propriétés de ces vecteurs. L'étude consiste donc à déterminer tous les nombres caractéristiques des vecteurs de parité. L'étude des propriétés des **nombres caractéristiques** nous permet de déterminer un certain nombre des équations ou des conditions que déterminent si un vecteur de parité est réalisable ou non réalisable. Noter qu'un vecteur de parité V de **longueur infinie** est dit **réalisable** s'il existe une suite de Collatz de longueur infinie et de premier terme un entier naturel non nul bien déterminé qui admet V comme vecteur de parité donc le cas contraire on dit que V est **non réalisable**.

La fonction de Collatz est définie comme suit pour tout entier naturel non nul :

$$\forall N \in \mathbb{N}^*, T(N) = \begin{cases} \dfrac{N}{2} & \text{si } N \equiv 0 \bmod(2) \\ \dfrac{3N+1}{2} & \text{si } N \equiv 1 \bmod(2) \end{cases}$$

**Notation 1.1.** La suite de Collatz de premier terme N et de longueur n sera notée comme suit :
$$Sy_n(N) = (N, T(N), T^2(N), \ldots, T^{n-1}(N))$$

**Notation 1.2.** Le vecteur de parité de la suite de Collatz de premier terme N et de longueur n est noté comme suit:
$$v_n(N) = (\mathbf{i}(N), \mathbf{i}(T(N)), \ldots, \mathbf{i}(T^{n-1}(N)))$$

Sachant que **i** est une application de $\mathbb{N}$ dans {0,1} définie comme suit:
$$\forall\, N \in \mathbb{N}, \mathbf{i}(N) = \begin{cases} 0 & \text{si } N \equiv 0 \mod(2) \\ 1 & \text{si } N \equiv 1 \mod(2) \end{cases}$$

Dans le cas ou le vecteur de parité de la suite de Collatz de premier terme N est de longueur infinie, il sera noté comme suit :
$$v_\infty(N) = (\mathbf{i}(N), \mathbf{i}(T(N)), \mathbf{i}(T^2(N)), \ldots, \mathbf{i}(T^n(N)), \ldots)$$

**Définition 1.2.** Soit **V** un vecteur binaire de **longueur infinie** alors on peut distinguer deux cas différents :

- On dit que **V** est **réalisable ou encore convertible** en suite de Collatz s'il existe un entier naturel non nul N tel que
$$v_\infty(N) = \mathbf{V} \qquad (\mathbf{1.1})$$

- On dit qui **V** est **non réalisable ou non convertible** en suite de Collatz si on a :
$$\forall\, N \in \mathbb{N}^* \quad v_\infty(N) \neq \mathbf{V} \qquad (\mathbf{1.2})$$

**Exemple 1.1.** On considère les deux vecteurs de parité de longueurs infinies comme suit :

Le premier vecteur est noté $V_1$ et il est défini comme suit :
$$V_1 = (1,1,1,0,1,0,0,1,0,0,0,1,0,1,0,1,0,\ldots,1,0,\ldots)$$

Ce vecteur binaire est considéré comme réalisable ou convertissable en suite de Collatz en effet la suite de Collatz de premier terme 7 et de longueur infinie a pour vecteur de parité le vecteur $V_1$ c'est à dire que :
$$v_\infty(7) = V_1$$

Le vecteur binaire noté $V_2$, il est comme suit :
$$V_2 = (1,0,0,1,0,0,1,0,0,1,0,0,1,0,0,\ldots,1,0,0\ldots)$$

Ce vecteur est considéré comme non convertible en suite de Collatz en effet si on suppose qu'il existe un entier naturel non nul x tel que $v_\infty(x) = V_2$ alors nécessairement x vérifie l'équation suivante :

$$T^3(x) = x \in \mathbb{N}$$

Equivaut à :

$$T^3(x) = \frac{3}{2^3}x + \frac{1}{2^3} = x \Rightarrow x = \frac{1}{5} \notin \mathbb{N}^*$$

On conclut qu'il n'existe pas aucun **entier naturel non nul** x tel que la suite de Collatz de premier terme x et de **longueur infinie** admettant $V_2$ comme vecteur de parité autrement :

$$\forall\, x \in \mathbb{N}^*, v_\infty(x) \neq V_2$$

**Notations 1.3.** On désigne par :

$\mathbb{B}$: L'ensemble constitué par tous les vecteurs binaires de **longueurs infinies**,

$\mathbb{B}_R$ : La partie de $\mathbb{B}$ qui contient uniquement les vecteurs de parité infinis **réalisable ou convertible** en suites de Collatz autrement :

$$\mathbb{B}_R = \{v_\infty(N) \mid N \in \mathbb{N}^*\} \quad (1.3)$$

$\mathbb{B}_{NR}$ : La partie de $\mathbb{B}$ qui contient tous les vecteurs de parité infinis **non réalisables ou non convertibles** en suites de Collatz c'est le complément de $\mathbb{B}_R$ c'est à dire que $\mathbb{B}_{NR} = \overline{\mathbb{B}_R}$ alors on peut écrire :

$$\mathbb{B} = \mathbb{B}_R \cup \mathbb{B}_{NR} \quad (1.4)$$

$\mathbb{B}_{(0,1)}$: L'ensemble qui renferme tous les vecteurs binaires infinis tel que chaque vecteur est constitué d'**une partie finie** que se comporte comme un vecteur de parité quelconque de longueur fini suivi d'une partie périodique de longueur infinie dont le cycle est (0,1) autrement si on désigne par V un vecteur binaire infini appartenant à $\mathbb{B}_{(0,1)}$ donc V sera de la forme suivante :

$$V = (\underbrace{e_1, e_2, \ldots, e_n}_{\text{Portion de longueur finie}}, \underbrace{0,1,0,1,0,1, \ldots, 0,1, \ldots}_{\text{Partie périodique de longueur infinie}}) \quad (1.5)$$

$\forall i \in [\![1, n]\!], e_i \in \{0,1\}$ et $n \in \mathbb{N}^*$

**Enoncée de la conjecture de Collatz**

Toutes les suites infinies qui sont générées par la fonction de Collatz et dont leurs premiers termes sont des entiers naturels non nuls admettant des vecteurs de parité appartenant à l'ensemble $\mathbb{B}_{(0,1)}$ autrement dit:

$$\forall\, N \in \mathbb{N}^*, v_\infty(N) \subset \mathbb{B}_{(0,1)} \quad (1.6)$$

**Selon la conjecture de Collatz**, ceci signifie que:

$$\mathbb{B}_R \subset \mathbb{B}_{(0,1)} \quad (1.7)$$

## 2. Généralité sur la caractérisation des vecteurs de parité

Si on part maintenant d'un vecteur de parité quelconque noté v, de longueur n avec n un entier naturel non nul. On sait qu'un tell vecteur de parité est un élément de $\{0,1\}^n$ donc il peut être représenté comme suit :

$$v = (e_1, e_2, \ldots, e_n) \qquad (2.1)$$

Tel que $\forall\, i \in [\![1, n]\!], e_i \in \{0,1\}$

**Notation 2.1.** Soit v un vecteur de parité quelconque de longueur finie n, on sait qu'elle existe une infinité des suites de Collatz de même longueur n admettant v comme vecteur de parité alors on désigne par $\mathbf{N_0(v)}$ le plus entier naturel non nul tel que la suite de Collatz de premier terme $N_0(v)$ et de longueur n admet v comme vecteur de parité. Autrement dit :

$$\begin{cases} v_n(N_0(v)) = v \\ \forall\, x \in \mathbb{N}^* \text{ tel que } v_n(x) = v, \text{ on a } x \geq N_0(v) \end{cases} \qquad (2.2)$$

Noter que :

$$1 \leq N_0(v) \leq 2^{n(v)}$$

**Exemple 2.1.** On considère le vecteur de parité suivant :

$$v = (1,0,1,1,1,0)$$

9 est le plus petit entier naturel non nul tel que :

$$v_6(9) = v$$

On déduit que :

$$N_0(v) = 9$$

**Notation 2.2.** Soit v un vecteur de parité de longueur finie alors pour tout entier naturel j, on désigne par $N_j(v)$ l'entier naturel non nul définie par l'expression suivante :

$$N_j(v) = N_0(v) + 2^{n(v)} j \qquad (2.3)$$

Il clair que :

$$N_{j+1}(v) - N_j(v) = 2^{n(v)}$$

**Remarque 2.1.** La suite de Collatz de premier terme $N_j(v)$ et de longueur n(v) admet v comme vecteur de parité.

**Notation 2.3.** On définit l'ensemble constitué par tous les entiers naturels non nuls tel que chaque élément représente le premier terme d'une suite de Collatz dont le vecteur de parité est le vecteur v.

$$\mathbb{H}(v) = \left\{N_j(v) \in \mathbb{N}^* \mid j \in \mathbb{N}, v_{n(v)}\left(N_j(v)\right) = v\right\} \quad (2.4)$$

Ou encore :

$$\mathbb{H}(v) = \left\{N_0(v), N_1(v), \dots, N_j(v), \dots\right\}$$

Toute les suites de Collatz dont les premiers termes sont des éléments de $\mathbb{H}(v)$ et de longueur $n(v)$ admet $v$ comme vecteur de parité.

**Notation 2.4.** Soit $v$ un vecteur de parité de longueur finie $n$. On désigne par $S_0(v)$ la suite de Collatz de longueur $n(v)$, qui a pour vecteur de parité le vecteur $v$ et de premier terme $N_0(v)$ donc la notation $S_0(v)$ ne qu' une autre notation de la suite $Sy_{n(v)}(N_0(v))$ c'est à dire que :

$$S_0(v) = Sy_{n(v)}\left(N_0(v)\right) \quad (2.5)$$

La suite de premier terme $N_j(v)$ et de longueur $n(v)$ est notée $S_j(v)$, elle admet $v$ comme vecteur de parité.

$$S_j(v) = Sy_{n(v)}(N_j(v))$$

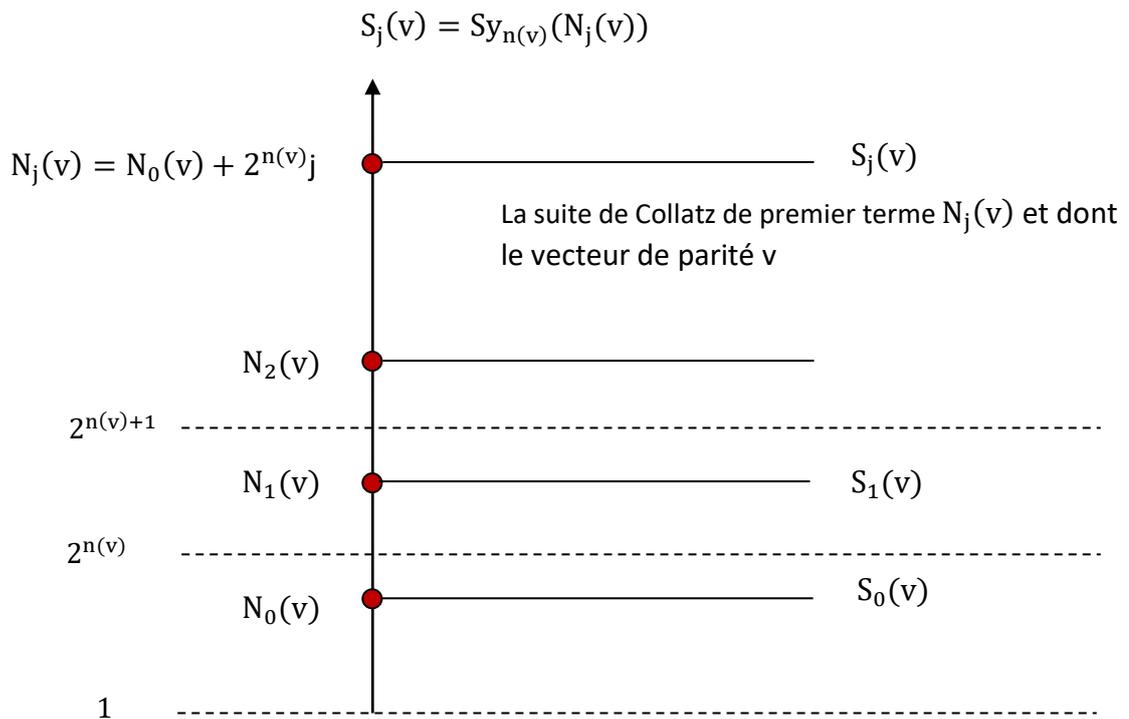

**FiG 2 :** Répartition des suites $S_j(v)$

**Définition 2.1.** Soit $v$ un élément de $\{0,1\}^n$ avec $n$ un entier naturel non nul, **un nombre caractéristique de v** est tout nombre (entier ou non entier) qui intervient dans la caractérisation du vecteur $v$ ou encore on peut le définir comme tout nombre qui nous permet de déterminer une ou plusieurs propriétés du $v$. Les nombres caractéristiques sont notés généralement sous la forme d'une lettre suivi par la lettre $v$ placée entre deux parenthèses par exemple $x(v), y(v), a(v),\dots$

**Exemple 2.2.** Soient n un entier naturel non nul, m un entier naturel quelconque et on considère le vecteur binaire v comme suit :

$$v = (e_1, e_2, \ldots, e_n)$$

Tel que $\forall\ i \in [\![1, n]\!], e_i \in \{0,1\}$

On suppose que :

$$\sum_{i=1}^{n} e_i = m \tag{2.6}$$

La longueur n de vecteur binaire v et le nombre m figurent parmi les nombres caractéristiques du vecteur v et ils seront notés **n(v) et m(v).**

**Définition 2.2.** Soit v un vecteur de parité de longueur finie. Si on désigne par $x(v), y(v), z(v), ..$ un ensemble des nombres caractéristiques quelconques de v alors par définition une équation caractéristique de v est toute équation décrivant une relation entre deux ou plusieurs nombres caractéristiques de v autrement dit si f est une fonction définie sur $\mathbb{R}$ dans $\mathbb{R}$ alors **une équation caractéristique de v** s'écrit sous la forme suivante :

$$f(x(v), y(v), z(v), \ldots) = 0$$

**Remarque 2.2.** Chaque vecteur v peut être caractérisé par plusieurs nombres caractéristiques (nombres entiers ou non entiers). Les nombres étudiés dans cet article sont notés comme suit :

$$n(v), m(v), P(v), a(v), b(v), X(v), Y(v), c(v), \ldots$$

Par exemple, les trois nombres caractéristiques **n(v), m(v)** et **P(v)** sont les nombres caractéristiques fondamentaux ou de base et ils sont ceux qui apparaissent dans l'équation suivante :

$$\mathbf{T^{n(v)}(N) = \frac{3^{m(v)}}{2^{n(v)}} N + \frac{P(v)}{2^{n(v)}}} \tag{2.8}$$

Les deux nombres **a(v)** et **b(v)** sont deux entiers naturels non nus qui vérifient l'équation suivante appelée équation caractéristique du vecteur de parité v :

$$\mathbf{3^{m(v)} a(v) + 1 = 2^{n(v)} b(v)} \tag{2.9}$$

**Définition 2.3.** On sait que pour un vecteur de parité fini v, elle existe une infinité des suites admettant v comme vecteur de parité, les premiers termes des ces suites sont notés $N_0(v), N_1(v), ..$ on suppose qu'il existe un terme de cet ensemble qui peut être exprimé en fonction d'un ou plusieurs nombres caractéristiques de v donc ce nombre

$$N_j(v) = f(\text{nombres caractéristiques})$$

Dans ce cas, le terme $N_j(v)$ est appelle point particulier ou premier terme caractéristique du vecteur de parité v.

**Exemple 2.3.** Si l'équation caractéristique précédente est multipliée par P(v) on obtient:

$$3^{m(v)}a(v)P(v) + P(v) = 2^{n(v)}b(v)P(v)$$

La division de cette dernière par $2^{n(v)}$ nous permet d'écrire :

$$\frac{3^{m(v)}}{2^{n(v)}}a(v)P(v) + \frac{P(v)}{2^{n(v)}} = b(v)P(v)$$

Cette dernière signifie que la suite de Collatz de premier terme **X(v)=a(v)P(v)** et de longueur n(v) admet v comme vecteur de parité et dans ce cas on peut écrire:

$$\frac{3^{m(v)}}{2^{n(v)}}X(v) + \frac{P(v)}{2^{n(v)}} = Y(v)$$

Avec **Y(v) = b(v)P(v)**

**Notation 2.5.** Un vecteur de parité de longueur infinie peut être représenté comme suit :

$$V = (e_1, e_2, \ldots, e_j, e_{j+1}, \ldots) \; \forall \; i \in \mathbb{N}^*, e_i \in \{0,1\}$$

Pour la caractérisation de V, on désigne par $R_j(V)$ la partie de V qui contient les j premiers termes de V c'est à dire que :

$$R_j(V) = (e_1, e_2, \ldots, e_j) \text{ avec } j \in \mathbb{N}^* \tag{2.10}$$

**Corollaire 2.1.** Soit V un vecteur de parité de longueur infinie donc il est évident que :

$$\lim_{j \to +\infty} R_j(V) = V \tag{2.11}$$

**Notation 2.6.** L'ensemble ordonné constitué par toutes les parties de V est noté comme suit :

$$\mathbf{D(V) = (R_1(V), R_2(V), R_3(V), \ldots, R_j(V), \ldots)} \tag{2.12}$$

**Définition 2.4.** Si on désigne par $x\big(R_j(V)\big)$ un nombre caractéristique quelconque de $R_j(V)$ alors ce même nombre représente aussi **un nombre caractéristique d'ordre j du vecteur V**. Dans ce cas le nombre caractéristique d'ordre j de V est noté comme suit :

$$\mathbf{x\big(R_j(V)\big) = x_j(V)} \tag{2.13}$$

**Exemple 2.4.** Si on considère le vecteur V suivant supposé de **longueur infinie** :

$$V = (1,1,0,1,0,0,1,0,1,1,1,0,0, \ldots)$$

La partie de V constituée par ses 6 premiers termes est comme suit :

$$\mathbf{R_6(V) = (1, 1, 0, 1, 0, 0)}$$

Le nombre m(v) désigne le nombre de termes non nul dans un vecteur de parité quelconque alors on peut écrire :

$$\mathbf{m}\big(\mathbf{R_6}(\mathbf{V})\big) = \mathbf{3}$$

Ce dernier est un nombre caractéristique d'ordre 6 du vecteur V considéré alors on écrit :

$$\mathbf{m_6}(\mathbf{V}) = \mathbf{3}$$

**Définition 2.5.** Soit V un vecteur binaire de longueur infinie et on désigne par $x_j(V)$ un nombre caractéristique d'ordre j de vecteur V avec j un entier naturel non nul quelconque. Alors la suite constituée par tous les nombres caractéristiques d'ordre j allant de 1 à l'infini constitue une suite caractéristique du vecteur V, elle notée comme suit :

$$\mathbf{L_x(V)} = \big(\mathbf{x_1(V), x_2(V), x_3(V), \ldots, x_j(V), \ldots}\big) \qquad (2.14)$$

**Définition 2.6.** Les nombres caractéristiques absolus (ou tout simplement les nombres caractéristiques) d'un vecteur de parité de longueur infinie sont les limites des nombres caractéristiques d'ordre finie autrement dit si on désigne par $x_j(V)$ un nombre caractéristique d'ordre j de vecteur V alors le nombre caractéristique absolus de V noté x(V) est défini comme suit :

$$\mathbf{x(V)} = \lim_{\mathbf{j} \to +\infty} \mathbf{x_j(V)} \qquad (2.15)$$

**Définition 2.7.** Soient V un vecteur binaire de **longueur infinie** et j un entier naturel non nul alors toute équation caractéristique relative au vecteur de parité $R_j(V)$ représente une **équation caractéristique d'ordre j** du vecteur binaire infini V.

**Exemple 2.5.** Si un vecteur de parité de longueur finie v est caractérisé par l'équation suivante :

$$3^{m(v)}a(v) + 1 = 2^{n(v)}b(v)$$

Soit **V** un vecteur binaire de longueur infinie et de plus on suppose que v vérifie la relation suivante:

$$v = R_j(\mathbf{V})$$

Alors tous les nombre caractéristique de v sont des nombres caractéristiques d'ordre j de **V**, ce qui nous permet d'écrire :

$$a(v) = a\big(R_j(\mathbf{V})\big) = a_j(\mathbf{V}) \qquad (2.16)$$

L'équation ci-dessous est une **équation caractéristique d'ordre j** de V, elle s'écrit comme suit :

$$3^{m_j(V)} a_j(V) + 1 = 2^{n_j(V)} b_j(V) \tag{2.17}$$

**Notation 2.7.** On sait que $\mathbb{H}(R_j(V))$ désigne l'ensemble qui contient tous les premiers termes de toutes les suites de Collatz qui admettant le même vecteur de parité $R_j(V)$ cet ensemble constitue l'ensemble d'ordre j de vecteur binaire infinie V et on peut le noter comme suit :

$$\mathbb{H}\left(R_j(V)\right) = \mathbb{H}_j(V) \tag{2.18}$$

Rappeler que :

$$\mathbb{H}\left(R_j(V)\right) = \left\{N_0\left(R_j(V)\right), N_1\left(R_j(V)\right), \ldots, N_k\left(R_j(V)\right), \ldots\right\}$$

Le premier terme de la $k^{\text{éme}}$ suite de Collatz dont le vecteur de parité est le vecteur $R_j(V)$ est noté comme suit:

$$N_k\left(R_j(V)\right) = N_{k,j}(V) \tag{2.19}$$

Le plus petit entier naturel non nul qui constitue le premier terme de la suite de Collatz de vecteur de parité $R_j(V)$ est noté comme suit:

$$N_0\left(R_j(V)\right) = N_{0,j}(V) \tag{2.20}$$

**Notation 2.8.** D'après les deux équations (2.18) et (2.19), on peut écrire :

$$\mathbb{H}_j(V) = \left\{N_{0,j}(V), N_{1,j}(V), \ldots, N_{k,j}(V), \ldots\right\} \tag{2.21}$$

**Notation 2.9.** La suite dont les termes sont les valeurs de $N_{0,j}(V)$ pour $j \in \mathbb{N}^*$ est notée comme suit :

$$h_0(V) = \left(N_{0,1}(V), N_{0,2}(V), N_{0,3}(V), \ldots, N_{0,j}(V), \ldots\right) \tag{2.22}$$

Noter que :

$$\forall j \in \mathbb{N}^*, \quad 1 \leq N_{0,j}(V) \leq 2^{n_j(V)} \tag{2.23}$$

**Notation 2.10.** Soit V un vecteur de parité quelconque de longueur infinie alors pour tout entier naturel non nul j, On désigne par $S_{0,j}(V)$, la suite de Collatz de premier terme $N_{0,j}(V)$ et de longueur $n_j(V)$. Noter que $R_j(V)$ est le vecteur de parité de cette suite.

**Exemple 2.6.** On considère le vecteur binaire V de **longueur infinie**, on suppose que ses 14 premiers termes sont comme suit :

$$V = (1,1,0,1,0,0,1,1,0,1,0,0,1,0,\ldots)$$

On cherche à déterminer les 8 premiers termes des suites $S_{0,j}(V)$ pour des valeurs de j allants de 1 à 8. Rappeler que $N_{0,j}(V)$ est le plus petit entier naturel non nul qui représente le premier terme d'une suite de Collatz de longueur j ($n_j(V) = j$) qui admet $R_j(V)$ comme vecteur de parité.

Par exemple 11 est le plus entier naturel qui représente le premier terme d'une suite de Collatz de longueur 4 et de vecteur de parité $R_4(V) = (1,1,0,1)$ alors que 139 est le plus entier naturel qui représente le premier terme d'une suite de Collatz de longueur 8 et de vecteur de parité $R_8(V) = (1,1,0,1,0,0,1,1)$ donc on peut écrire :

$$N_{0,4}(V) = 11, N_{0,8}(V) = 139$$

Les expressions de $R_j(V)$ pour j allant de 1 à 8 sont comme suit :

$$R_1(V) = 1, R_2(V) = (1,1), R_3(V) = (1,1,0), R_4(V) = (1,1,0,1), R_5(V) = (1,1,0,1,0)$$
$$R_6(V) = (1,1,0,1,0,0), R_7(V) = (1,1,0,1,0,0,1), R_8(V) = (1,1,0,1,0,0,1,1)$$

Le tableau ci-dessous contient les valeurs de $N_{0,j}(V)$ pour j allant de 1 à 8 qui sont comme suit :

**Tableau 1. Les suites $S_{0,j}(V)$**

| 1 | 1 | 0 | 1 | 0 | 0 | 1 | 1 | 0 | 1 |
|---|---|---|---|---|---|---|---|---|---|
| 1 | | | | | | | | | |
| 3 | 5 | | | | | | | | |
| 3 | 5 | 8 | | | | | | | |
| 11 | 17 | 26 | 13 | | | | | | |
| 11 | 17 | 26 | 13 | 20 | | | | | |
| 11 | 17 | 26 | 13 | 20 | 10 | | | | |
| 11 | 17 | 26 | 13 | 20 | 10 | 5 | | | |
| 139 | 209 | 314 | 157 | 236 | 118 | 59 | 89 | | |

On déduit les différentes valeurs de $N_{0,j}(V)$ pour $j \in [\![1,8]\!]$ comme suit :

$N_{0,1}(V) = 1$

$N_{0,2}(V) = N_{0,3}(V) = 3$

$N_{0,4}(V) = N_{0,5}(V) = N_{0,6}(V) = N_{0,7}(V) = 11$

$N_{0,8}(V) = 139$

**Remarque 2.3.** On suppose que :

$$\lim_{j \to +\infty} N_{0,j}(V) = +\infty$$

Dans ce cas on peut conclure qu'elle n'existe pas aucune suite de Collatz de **longueur infinie** qui admet V comme vecteur de parité. Le vecteur binaire de longueur infinie V est considéré comme **non réalisable ou non convertible** en suite d Collatz dans ce cas.

**Proposition 2.1.** Soit V un vecteur binaire de **longueur infinie** donc peut distinguer deux cas possibles pour l'existence ou le non existence d'une suite de Collatz de longueur infinie qui admet V comme vecteur de parité. Si :

$$\lim_{j \to +\infty} N_{0,j}(V) = +\infty \qquad (2.24)$$

Dans ce cas le vecteur binaire V est non réalisable et elle n'existe pas aucune suite de Collatz de premier terme bien déterminé qui admet V comme vecteur de parité.

$$\lim_{j \to +\infty} N_{0,j}(V) \text{ est finie} \qquad (2.25)$$

Dans ce cas, le vecteur binaire V est réalisable et elle existe une suite de Collatz de longueur infinie et de premier terme un entier naturel non nul bien déterminé qui correspond a $N_{0,\infty}(V)$ tel que le vecteur de parité de cette suite est le vecteur binaire infini V.

## 3. Suites associées aux suites de Collatz

**Définition 3.1.** Soit m un entier naturel non nul bien déterminé (m ∈ ℕ*), on définit les deux suites $(a_{m,n})_{n \geq 1}$ et $(b_{m,n})_{n \geq 1}$ comme suit :

Pour n=1

$$\begin{cases} a_{m,1} = 1 \\ b_{m,1} = \dfrac{3^m + 1}{2} \end{cases} \qquad (3.1)$$

∀n ≥ 2

-Si $b_{m,n-1}$ est pair

$$\begin{cases} a_{m,n} = a_{m,n-1} \\ b_{m,n} = \dfrac{b_{m,n-1}}{2} \end{cases} \qquad (3.2)$$

-Si $b_{m,n-1}$ est impair

$$\begin{cases} a_{m,n} = a_{m,n-1} + 2^{n-1} \\ b_{m,n} = \dfrac{b_{m,n-1} + 3^m}{2} \end{cases} \qquad (3.3)$$

**Corollaire 3.1.** Soient m un entier naturel non nul quelconque alors pour tout entier naturel non nul n la relation entre les deux termes $a_{m,n}$ et $b_{m,n}$ de deux suites $(a_{m,n})_{n \geq 1}$ et $(b_{m,n})_{n \geq 1}$ est décrite par l'équation suivante :

$$3^m a_{m,n} + 1 = 2^n b_{m,n} \qquad (3.4)$$

**Démonstration.** On montre par récurrence que cette propriété est vraie pour tout entier naturel non nul n pour une valeur de m donnée.

Pour n =1, on sait que $3^m + 1$ est un entier naturel non nul pair donc il existe un entier naturel non nul qu'on peut noter $b_{m,1}$ tel que:

$$3^m + 1 = 2b_{m,1}$$

Donc quelque soit la valeur non nul de m la propriété est vraie pour n=1 et dans ce cas on peut écrire :

$$\begin{cases} a_{m,1} = 1 \\ b_{m,1} = \dfrac{3^m + 1}{2} \end{cases}$$

On suppose que la propriété est vraie pour tout entier naturel k allant de 2 à n et on montre que la propriété est vraie pour k=n+1. On peut écrire alors:

$$3^m a_{m,n} + 1 = 2^n b_{m,n}$$

On peut distinguer deux cas :

-Si $b_{m,n}$ est pair donc dans ce cas on peut écrire :

$$b_{m,n+1} = \dfrac{b_{m,n}}{2}$$

Ce qui implique que :

$$3^m a_{m,n} + 1 = 2^{n+1} b_{m,n+1}$$

On conclut que dans ce cas la propriété est vraie pour n+1 et on a :

$$\begin{cases} a_{m,n+1} = a_{m,n} \\ b_{m,n+1} = \dfrac{b_{m,n}}{2} \end{cases}$$

-Si $b_{m,n}$ est impair donc dans ce cas, on peut écrire en ajoutant $2^n a_{m,n}$ de deux cotés dans l'équation comme suit:

$$3^m a_{m,n} + 2^n 3^m + 1 = 2^n b_{m,n} + 2^n 3^m$$

Equivalente à :

$$3^m (a_{m,n} + 2^n) + 1 = 2^n (b_{m,n} + 3^m)$$

Comme $(b_{m,n} + 3^m)$ est pair donc il existe un entier naturel non nul $b_{m,n+1}$ tel que :

$$b_{m,n} + 3^m = 2 b_{m,n+1}$$

Ceci nous permet d'écrire :

$$3^m (a_{m,n} + 2^n) + 1 = 2^{n+1} b_{m,n+1}$$

On peut conclure que la propriété est vraie aussi dans ce cas pour k=n+1 et on a :

$$\begin{cases} a_{m,n+1} = a_{m,n} + 2^n \\ b_{m,n+1} = \dfrac{b_{m,n} + 3^m}{2} \end{cases}$$

Par conclusion la propriété est vraie pour tout entier naturel non nul n.

**Corollaire 3.2.** Soit m et n deux entiers naturels non nuls, on sait que la relation entre $a_{m,n}$ et $b_{m,n}$ se traduit par l'équation suivante :

$$3^m a_{m,n} + 1 = 2^n b_{m,n}$$

Alors les deux entiers naturels $a_{m,n}$ et $b_{m,n}$ vérifient les conditions suivantes :

$$a_{m,n} < 2^n \text{ et } b_{m,n} < 3^m \qquad (3.5)$$

**Démonstration** On sait que :

Pour i=1

$$a_{m,i} = 1$$

$\forall\, i \geq 2$

$$\begin{cases} a_{m,i} = a_{m,i-1} & \text{si } b_{m,i-1} \text{ est pair} \\ a_{m,i} = a_{m,i-1} + 2^{i-1} & \text{si } b_{m,i-1} \text{ est impair} \end{cases}$$

On déduit l'expression de $a_{m,n}$ comme suit :

$$a_{m,n} = 1 + 2\mathbf{i}(b_{m,1}) + 2^2\mathbf{i}(b_{m,2}) + 2^3\mathbf{i}(b_{m,3}) + \cdots + 2^{n-1}\mathbf{i}(b_{m,n-1})$$

Rappelons que :

$$\begin{cases} \mathbf{i}(b_{m,j}) = 1 & \text{si } b_{m,j} \text{ est impair} \\ \mathbf{i}(b_{m,j}) = 0 & \text{si } b_{m,j} \text{ est pair} \end{cases}$$

**La valeur maximale** de $a_{m,n}$ s'écrit comme suit :

$$a_{m,n} = 1 + 2 + 2^2 + 2^3 + \cdots + 2^{n-1}$$

$$= 2^n - 1$$

La valeur maximale de $a_{m,n}$ est strictement inférieure a $2^n$.

**Exemple 3.1.** On considère le vecteur de parité v ci-dessous de longueur n(v)=10 :

| 1 | 0 | 1 | 1 | 0 | 1 | 0 | 1 | 1 | 1 |
|---|---|---|---|---|---|---|---|---|---|

Le vecteur de parité considéré est caractérisé par les nombres suivants :

$$\begin{cases} n(v) = 10 \\ m(v) = 7 \\ P(v) = 5645 \end{cases}$$

On cherche à déterminer deux entiers naturels non nul $a_{7,10}$ et $b_{7,10}$ tel que :

$$3^{m(v)} a_{7,10} + 1 = 2^{n(v)} b_{7,10}$$

Ce qui équivaut à :

$$3^7 a_{7,10} + 1 = 2^{10} b_{7,10}$$

On adopte la méthode suivante :

On sait que :

$$3^7 + 1 = 2 \times 1094 = 2^2 \times 547$$

Ce qui signifie que :

$$a_{7,1} = a_{7,2} = 1, \qquad b_{7,1} = 1094, \qquad b_{7,2} = 547$$

On ajoute $3^7 \text{x} 2^2$ on obtient :

$$3^7(1 + 2^2) + 1 = 2^2 \text{ x } (547 + 3^7) = 2^3 \text{x } 1367$$

On déduit alors :

$$a_{7,3} = 3, \qquad b_{7,3} = 1367$$

Puis on ajoute $3^7 \text{x} 2^3$ de deux cotés ce qui nous permet d'obtenir :

$$3^7(1 + 2^2 + 2^3) + 1 = 2^3 \text{x}( 1367 + 3^7) = 2^4 \text{x} 1777$$

Ce qui nous permet d'écrire :

$$a_{7,4} = 13, \qquad b_{7,4} = 1777$$

On ajoute alors $3^7 \text{x} 2^4$, on obtient :

$$3^7(1 + 2^2 + 2^3 + 2^4) + 1 = 2^4 \text{x}(1777 + 3^7) = 2^5 \text{x} 1982 = 2^6 \text{x} 991$$

On déduit que :

$$a_{7,5} = a_{7,6} = 29, b_{7,5} = 1982, b_{7,6} = 991$$

On ajoute alors $3^7 \text{x} 2^6$, on obtient :

$$3^7(1 + 2^2 + 2^3 + 2^4 + 2^6) + 1 = 2^6 \text{x}(991 + 3^7) = 2^7 \text{x} 1589$$

Ceci nous donne :

$$a_{7,7} = 93, \qquad b_{7,7} = 1589$$

On continu le calcul, on obtient :

$$3^7(1 + 2^2 + 2^3 + 2^4 + 2^6 + 2^7) + 1 = 2^8 \text{ x} 1888 = 2^9 \text{ x} 944 = 2^{10} \text{ x} 472$$

Il en résulte alors :

$$a_{7,8} = a_{7,9} = a_{7,10} = 221; \; b_{7,8} = 1888, b_{7,9} = 944, b_{7,10} = 472$$

On déduit que :

$$\begin{cases} a(v) = 221 \\ b(v) = 472 \end{cases}$$

On déduit que :

$$\mathbf{i}(b_{7,2}) = \mathbf{i}(b_{7,3}) = \mathbf{i}(b_{7,4}) = \mathbf{i}(b_{7,6}) = \mathbf{i}(b_{7,7}) = \mathbf{1}$$

On sait que l'expression de a(v) est comme suit :

$$a(v) = 1 + 2^2 + 2^3 + 2^4 + 2^6 + 2^7$$

**Corollaire 3.3.** Soient m, n deux entiers naturels non nuls quelconques et j un entier relatif non nul alors ils existent une infinité des couples de deux entiers relatifs non nuls qu'on peut noter $a_{m,n,j}$ et $b_{m,n,j}$ tel que :

$$3^m a_{m,n,j} + 1 = 2^n b_{m,n,j} \tag{3.6}$$

Avec :

$$\begin{cases} a_{m,n,j} = a_{m,n} + 2^n j \\ b_{m,n,j} = b_{m,n} + 3^m j \end{cases} \quad (3.7)$$

$(m, n) \in \mathbb{N}^* \times \mathbb{N}^*$ et $j \in \mathbb{Z}^*$

**Démonstration**  On sait qu'ils existent deux entiers naturels non nuls tel que :

$$3^m a_{m,n} + 1 = 2^n b_{m,n}$$

Pour un entier relatif non nul quelconque j, on ajoute de deux cotes de l'équation précédente le nombre $2^n 3^m j$ ce qui nous permet d'obtenir :

$$3^m a_{m,n} + 2^n 3^m j + 1 = 2^n b_{m,n} + 2^n 3^m j$$

Ou encore :

$$3^m (a_{m,n} + 2^n j) + 1 = 2^n (b_{m,n} + 3^m j)$$

On pose :

$$\begin{cases} a_{m,n,j} = a_{m,n} + 2^n j \\ b_{m,n,j} = b_{m,n} + 3^m j \end{cases}$$

## 4. Caractérisation des vecteurs de parités de longueurs finies

**Notation 4.1.**  Un vecteur de parité est un vecteur binaire qu'on peut la représenter dans le cas général comme une combinaison quelconque des 0 et des 1 :

$$v = (e_1, e_2, \ldots, e_n) \quad (4.1)$$

Les éléments $e_i$ ne peuvent prendre que l'une de deux valeurs 0 ou bien 1 ceci se traduit par :

$$\forall\, i \in [\![1, n]\!] \quad e_i \in \{0,1\}$$

**Notation 4.2.**  Soit v un vecteur de parité quelconque de longueur n avec n un entier naturel non nul alors pour tout entier i allant de 1 à n la notation $e_i(v)$ designe le $i^{eme}$ terme de vecteur de parité v. Si il n y a pas d'ambigüité on peut écrire tout simplement $e_i$ pour cette même signification.

**Exemple 4.1.**  On considère le vecteur de parité ci-dessous :

$$v = (1,0,0,1,1,0,1)$$

Donc on peut écrire :

$$e_1(v) = 1, e_2(v) = 0, e_3(v) = 0, \ldots$$

**Notations 4.3.**  Le nombre de terme ou la longueur de v est noté $n(v)$. Dans notre cas :

$$n(v) = n \quad (4.2)$$

Le nombre de termes qui sont égaux à 1 est noté $m(v)$, il correspond à la somme suivante :

$$m(v) = \sum_{i=1}^{n} e_i \tag{4.3}$$

Il est évident que :

$$m(v) \leq n(v)$$

**Définition 4.1.** Soit v un vecteur de parité de longueur finie. On sait qu'elle existe une seule suite de Collatz de longueur n(v) et dont le premier terme est compris entre 1 et $2^{n(v)}$ et admet v comme vecteur de parité. On définit le nombre rationnel $r_0(v)$ comme suit :

$$r_0(v) = \frac{N_0(v)}{2^{n(v)}} \tag{4.4}$$

Il est évident que ce nombre compris entre 0 et 1.

$$0 < r_0(v) \leq 1 \tag{4.5}$$

**Définition 4.2.** On définit le nombre c(v) appellé nombre de Collatz relatif au vecteur binaire v par la relation suivante :

$$c(v) = 2^{n(v)} - 3^{m(v)} \tag{4.6}$$

**Notation 4.4.** On désigne par $\mathbb{B}$ l'ensemble constitué par tous les vecteurs binaires c'est à dire que :

$$\mathbb{B} = \{v \in \{0,1\}^n | n \in \mathbb{N}^* \} \tag{4.7}$$

**Définition 4.3**

**Vecteurs de type $b^-$**

Ces vecteurs sont caractérisés par des nombres c(v) strictement négatifs, ils correspondent aux cas $2^{n(v)} < 3^{m(v)}$ c'est à dire que:

$$c(v) < 0 \tag{4.8}$$

**Vecteurs de type $b^+$**

Ces vecteurs sont caractérisés par des nombres C(v) strictement positifs ils correspondent aux cas $2^{n(v)} > 3^{m(v)}$ :

$$c(v) > 0 \tag{4.9}$$

**Notation 4.5.** L'ensemble qui contient tous les vecteurs binaires dont leurs sont des entiers négatifs est noté comme suit :

$$\mathbb{B}^- = \{v \in \mathbb{B} \mid c(v) < 0\}$$

L'ensemble qui contient tous les vecteurs binaires dont leurs sont des entiers positifs est noté comme suit :

$$\mathbb{B}^+ = \{v \in \mathbb{B} \mid C(v) > 0\}$$

**Définition 4.4.** Soient n, m deux entiers naturels non nuls tel que m ≤ n et v un vecteur de parité de longueur finie n contient m termes qui sont égaux à 1 qu'on peut le représenter comme suit :

$$v = \left(0, \ldots 0, e_{j_1}, 0 \ldots, 0, e_{j_2}, 0, \ldots, \ldots 0, \ldots, 0, e_{j_i}, 0, \ldots, 0, \ldots, 0, e_{j_m}, 0, 0, \ldots, 0\right) \quad (4.10)$$

Tel que :

$$\forall\, i \in [\![1, m]\!], e_{j_i} = 1 \quad (4.11)$$

Autrement :

$$e_{j_1} = e_{j_2} = \cdots = e_{j_i} = \cdots = e_{j_m} = 1$$

On définit le nombre caractéristique relatif au vecteur v qu'on le note P(v) par la relation suivante :

$$P(v) = 3^{m-1}x2^{j_1-1} + 3^{m-2}x2^{j_2-1} + \cdots + 3x2^{j_{m-1}-1} + 2^{j_m-1} \quad (4.12)$$

**Exemple 4.2.** On considère le vecteur de parité suivant :

$$v = (1,1,0,1,0,0,1)$$

Les deux nombres caractéristiques n(v) et m(v) de v sont les suivants :

$$\begin{cases} n(v) = 7 \\ m(v) = 4 \end{cases}$$

La valeur du nombre P(v) est déterminée comme suit :

$$P(v) = 3^{4-1}x2^{1-1} + 3^{4-2}x2^{2-1} + 3^{4-3}x2^{4-1} + 3^{4-4}x2^{7-1} = 133$$

**Théorème 4.1.** Soit v un vecteur de parité de longueur n avec n un entier naturel non nul. On suppose que v peut être représenté comme suit :

$$v = (e_1, e_2, \ldots, e_n)$$

Tel que :

$$\forall\, i \in [\![1, n]\!]\ e_i \in \{0,1\}$$

Alors pour déterminer la valeur de P(v), on procède comme suit en utilisant la relation de récurrence ci-dessous :

pour j = 1

$$\begin{cases} P_1(v) = 0 & \text{si } e_1 = 0 \\ P_1(v) = 1 & \text{si } e_1 = 1 \end{cases} \quad (4.13)$$

$\forall\, j \in [\![2, n]\!]$

$$P_j(v) = \begin{cases} P_{j-1}(v) & \text{si } e_j = 0 \\ 3P_{j-1}(v) + 2^{j-1} & \text{si } e_j = 1 \end{cases} \quad (4.20)$$

Le nombre caractéristique P(v) correspond au nombre caractéristique d'ordre n du vecteur v qu'on a noté $P_n(v)$ obtenue par itérations successives donc on peut écrire :

$$P_n(v) = P(v) \qquad (4.21)$$

**Exemple 4.3.** On reprend le vecteur de parité de l'exemple précédent et on fait calculer son nombre caractéristique en utilisant la relation de récurrence ci-dessus, les résultats obtenus sont représentés sur le tableau suivant :

| $e_j$ | 1 | 1 | 0 | 1 | 0 | 0 | 1 |
|---|---|---|---|---|---|---|---|
| $P_j(v)$ | 1 | 5 | 5 | 23 | 23 | 23 | 133 |

Le nombre caractéristique $P(v)$ est comme suit :

$$P(v) = P_7(v) = 133$$

**Théorème 4.2.** Soit N un entier naturel non nul, on désigne par $Sy_n(N)$ la suite de Collatz de premier terme N et de longueur n et par v le vecteur de parité de la suite considérée.

$$v_n(N) = v$$

Alors l'expression de $T^n(N)$ en fonction de N et des nombres caractéristiques de v est la suivante :

$$T^{n(v)}(N) = \frac{3^{m(v)}}{2^{n(v)}} N + \frac{P(v)}{2^{n(v)}} \qquad (4.22)$$

Noter que : $n = n(v)$

**Exemple 4.4.** On considère la suite de Collatz de premier terme 11 et de longueur 7 représentée ci-dessous :

$$Sy_7(11) = (11, 17, 26, 13, 20, 10, 5)$$

Le vecteur de parité de cette suite est le suivant :

$$v_7(11) = (1,1,0,1,0,0,1)$$

Les trois nombres caractéristiques de ce vecteur sont les suivants :

$$\begin{cases} n(v) = 7 \\ m(v) = 4 \\ P(v) = 133 \end{cases}$$

La valeur de $T^7(11)$ est obtenu par application de la formule nous donne :

$$T^{n(v)}(11) = T^7(11) = \frac{3^4}{2^7} 11 + \frac{133}{2^7} = 8$$

**Notation 4.6.** Soient m et n deux entiers naturels tel que $n \neq 0$ et $m \leq n$, on désigne par $\Omega(n,m)$ l'ensemble des vecteurs de parité de même longueur n et qui possèdent le même nombre caractéristique m(v) autrement dit :

$$\Omega(n, m) = \{v \in B \mid n(v) = n \text{ et } m(v) = m\} \qquad (4.23)$$

**Théorème 4.3.** Soient m et n deux entiers naturels non nuls tel que n ≤ m donc les valeurs qui correspondent aux tous les nombres caractéristiques P(v) de tous les vecteurs de parité appartenant à $\Omega(n, m)$ sont majorées par le nombre suivant :

$$\max\{P(v) \mid v \in \Omega(n, m)\} = 2^{n(v)-m(v)}\left(3^{m(v)} - 2^{m(v)}\right) \qquad (4.24)$$

Autrement dit :

$$\forall\, v \in \Omega(n, m),\, P(v) \leq 2^{n(v)-m(v)}\left(3^{m(v)} - 2^{m(v)}\right)$$

De plus le vecteur de parité de $\Omega(n, m)$ qui possède la valeur de P(v) **le plus grand** dans cet ensemble est défini comme suit :

$$v = (0, 0, 0, \ldots, 0, 0, 1, \ldots, 1, 1) \qquad (4.25)$$

Une telle distribution qu'on peut le décrire comme suit :

$$\forall\, i \in [\![1, n]\!]\ e_i = \begin{cases} 0 & \text{si } 1 \leq i \leq n - m \\ 1 & \text{si } n - m + 1 \leq i \leq n \end{cases}$$

**Démonstration.** On sait que le nombre caractéristique P(v) d'un vecteur de parité quelconque appartenant à $\Omega(n, m)$ peut être représentée dans le cas général sous la forme suivante :

$$P(v) = 3^{m-1}x2^{j_1 - 1} + 3^{m-2}x2^{j_2 - 1} + \cdots + 3x2^{j_{m-1} - 1} + 2^{j_m - 1}$$

$$\forall\, i \in [\![1, m]\!], e_{j_i} = 1$$

Pour obtenir la valeur maximale de P(v), on doit attribuer à chaque terme en exposant qu'on a notée $j_i$ sa valeur maximale ainsi la valeur maximale que peut prendre $j_m$ égale à n alors que la valeur maximale de $j_{m-1}$ est égal à n-1 et ainsi de suite donc on peut écrire :

$$\max(j_m) = n, \max(j_{m-1}) = n - 1, \ldots, \max(j_2) = n - m + 2, \max(j_1) = n - m + 1$$

L'expression du nombre caractéristique qui possède la valeur maximale dans (n,m) est la suivante :

$$P(v) = 3^{m-1}x2^{n-m} + 3^{m-2}x2^{n-m+1} + \cdots + 3x2^{n-2} + 2^{n-1}$$

Cette expression correspond à la somme dune suite géométrique de raison 3/2 de premier terme $2^{n-1}$ et de m termes

$$P(v) = 2^{n-1} \frac{\left(\frac{3}{2}\right)^m - 1}{\frac{3}{2} - 1} = 2^{n-m}(3^m - 2^m)$$

**Théorème 4.4.** Soient m et n deux entiers naturels non nuls tel que n ≥ m donc le minimum de tous les nombres caractéristiques P(v) de tous les vecteurs de parité appartenant à $\Omega(n, m)$ est donnée par l'expression suivante :

$$\min\{P(v) \mid v \in \Omega(n, m)\} = 3^{m(v)} - 2^{m(v)} \qquad (4.26)$$

Autrement dit :
$$\forall\, v \in \Omega(n, m), P(v) \geq 3^{m(v)} - 2^{m(v)} \quad (4.27)$$

De plus le vecteur de parité de $\Omega(n, m)$ qui possède la valeur $P(v)$ le plus petit dans cet ensemble est défini comme suit :
$$v = (1,1,1,\ldots 1,1,0,0,\ldots,0,0) \quad (4.28)$$

$$\forall\, i \in [\![1, n]\!] \ e_i = \begin{cases} 1 & \text{si } 1 \leq i \leq m \\ 0 & \text{si } m + 1 \leq i \leq n \end{cases}$$

**Démonstration.** On sait que le nombre $P(v)$ d'un vecteur de parité quelconque appartenant à $\Omega(n, m)$ peut être représenté dans le cas général sous la forme suivante :
$$P(v) = 3^{m-1} \times 2^{j_1 - 1} + 3^{m-2} \times 2^{j_2 - 1} + \cdots + 3 \times 2^{j_{m-1} - 1} + 2^{j_m - 1}$$

$$\forall\, i \in [\![1, m]\!], e_{j_i} = 1$$

Pour obtenir la valeur minimale de $P(v)$, on doit attribuer à chaque terme en exposant qu'on a notée $j_i$ sa valeur minimale ainsi la valeur minimale que peut prendre $j_1$ égale à 1 alors que la valeur minimale de $j_2$ est égal à 2 et ainsi de suite donc on peut écrire :

$$\min(j_1) = 1, \min(j_2) = 2, \ldots, \min(j_{m-1}) = m - 1, \min(j_m) = m$$

L'expression du nombre caractéristique qui possède la valeur minimale dans (n,m) est la suivante :
$$P(v) = 3^{m-1} + 3^{m-2} 2 + 3^{m-3} \times 2^2 + \cdots + 3 \times 2^{m-2} + 2^{m-1}$$

Cette expression correspond à la somme dune suite géométrique de raison 3/2 de premier terme $2^{m-1}$ et de m termes

$$P(v) = 2^{m-1} \frac{\left(\frac{3}{2}\right)^m - 1}{\frac{3}{2} - 1} = 3^m - 2^m$$

**Corollaire 4.1.** Tous les vecteur binaire v de longueurs finies vérifient l'inégalité suivante :
$$\frac{P(v)}{2^{n(v)} 3^{m(v)}} \leq \frac{1}{2^{m(v)}} - \frac{1}{3^{m(v)}} \quad (4.29)$$

**Proposition 4.1.** On considère un vecteur de parité v de longueur finie et on suppose que v peut être décomposé deux vecteurs $v_1$ et $v_2$ comme suit :

$$v = (\underbrace{e_1, e_2, \ldots, e_{n_1}}_{v_1}, \underbrace{e_{n_1+1}, e_{n_1+2}, \ldots, e_{n_1+n_2}}_{v_2})$$

Avec : $n(v_1) = n_1, n(v_2) = n_2,\ n(v) = n(v_1) + n(v_2)$

Alors l'expression du nombre caractéristique P(v) en fonction des nombres caractéristiques de $v_1$ et de $v_2$ s'écrit comme suit :

$$P(v) = 3^{m(v_2)}P(v_1) + 2^{n(v_1)}P(v_2) \qquad (4.30)$$

**Démonstration.** Soient $x_0, x_1$ et $x_2$ trois entiers naturels non nuls. On considère la suite de Collatz de premier terme $x_0$, de longueur $n(v_1)$ et de vecteur de parité $v_1$ et on suppose que :

$$T^{n(v_1)}(x_0) = x_1$$

Puis on considère la suite de Collatz de premier terme $x_1$, de longueur $n(v_2)$ et de vecteur de parité $v_2$ et on suppose que :

$$T^{n(v_2)}(x_1) = x_2$$

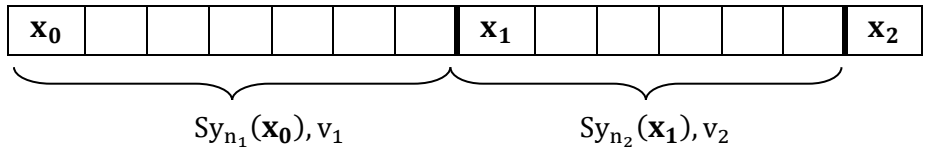

L'expression de $T^{n(v_1)}(x_0)$ en fonction de $x_0$ et de nombres caractéristiques de $v_1$ est comme suit :

$$T^{n(v_1)}(x_0) = \frac{3^{m(v_1)}}{2^{n(v_1)}} x_0 + \frac{P(v_1)}{2^{n(v_1)}}$$
$$= x_1$$

L'expression de $T^{n(v_1)}(x_1)$ en fonction de $x_1$ et de nombres caractéristiques de $v_2$ est comme suit :

$$T^{n(v_2)}(x_1) = \frac{3^{m(v_2)}}{2^{n(v_2)}} x_1 + \frac{P(v_2)}{2^{n(v_2)}}$$
$$= \frac{3^{m(v_2)}}{2^{n(v_2)}} \left( \frac{3^{m(v_1)}}{2^{n(v_1)}} x_0 + \frac{P(v_1)}{2^{n(v_1)}} \right) + \frac{P(v_2)}{2^{n(v_2)}}$$
$$= \frac{3^{m(v_2)+m(v_1)}}{2^{n(v_2)+n(v_1)}} x_0 + \frac{3^{m(v_2)}P(v_1) + 2^{n(v_1)}P(v_2)}{2^{n(v_1)+n(v_1)}}$$
$$= \frac{3^{m(v)}}{2^{n(v)}} x_0 + \frac{3^{m(v_2)}P(v_1) + 2^{n(v_1)}P(v_2)}{2^{n(v)}}$$
$$= \frac{3^{m(v)}}{2^{n(v)}} x_0 + \frac{P(v)}{2^{n(v)}}$$

On déduit l'expression de P(v) en fonction des nombres caractéristiques de $v_1$ et de $v_2$ qui s'écrit comme suit :

$$P(v) = 3^{m(v_2)}P(v_1) + 2^{n(v_1)}P(v_2)$$

**Corollaire 4.2.** Soit u un vecteur de parité quelconque de longueur finie n et qu'on peut le représenter comme suit :

$$u = (e_1, e_2, \ldots, e_n)$$

On considère un vecteur de parité périodique compose de k fois du vecteur de parité u comme suit :

$$v_k = (\underbrace{e_1, e_2, \ldots, e_n, \mathbf{e_1}, \mathbf{e_2}, \ldots, \mathbf{e_n}, e_1, e_2, \ldots, e_n, \ldots \ldots \ldots, \mathbf{e_1}, \mathbf{e_2}, \ldots, \mathbf{e_n}})$$
$$\text{K fois le vecteur u}$$

L'expression de $P(v_k)$ en fonction de $P(u)$ et des nombres caractéristiques de u est comme suit :

$$P(v_k) = \frac{3^{km(u)} - 2^{kn(u)}}{3^{m(u)} - 2^{n(u)}} P(u) \qquad (4.31)$$

**Exemple 4.5.** On considère le vecteur de parité suivant de longueur finie :

| 0 | 1 | 0 | 1 | 0 | 1 |  | 0 | 1 |

$$n(v) = 2m(v)$$

Le nombre caractéristique $P(v)$ de ce vecteur s'écrit comme suit :

$$P(v) = 2^{2m(v)} - 3^{m(v)}$$

**Définition 4.5.** Soit v un vecteur de parité quelconque de longueur finie donc ils existent deux entiers naturels non nuls qu'on peut les noter $a(v)$ et $b(v)$ tel que :

$$3^{m(v)} a(v) + 1 = 2^{n(v)} b(v) \qquad (4.32)$$

L'équation est appelée **l'équation caractéristique du vecteur de parité v.**

Les deux entiers $a(v)$ et $b(v)$ figurent parmi les nombres caractéristiques du vecteur v considéré. De plus ces deux nombres ne autres que les deux entiers naturels non nuls $a_{m,n}$ et $b_{m,n}$ donc on peut écrire :

$$\begin{cases} a(v) = a_{m,n} \\ b(v) = b_{m,n} \end{cases} \qquad (4.33)$$

**Définition 4.6.** Soit un vecteur de parité v de longueur finie, on définie les deux nombres caractéristiques $\alpha(v)$ et $\beta(v)$ qui sont deux entiers naturels non nuls tel que :

$$P(v) = 3^{m(v)} \alpha(v) + \beta(v) \mid \beta(v) < 3^{m(v)} \qquad (4.34)$$

De plus, on définie les deux nombres caractéristiques $A(v)$ et $B(v)$ qui sont deux entiers naturels non nuls tel que :

$$P(v) = 2^{n(v)} A(v) + B(v) \mid B(v) < 2^{n(v)} \qquad (4.35)$$

**Remarque 4.1.** Tous les vecteurs de parité de l'ensemble $\Omega(n, m)$ possédant la même équation caractéristiques.

**Définition 4.7.** Soit v un vecteur de parité quelconque, on définit la fonction notée $g_v$ de $\mathbb{N}^*$ dans $\mathbb{R}$ tel que :

$$g_v(N) = \frac{3^{m(v)}}{2^{n(v)}} N + \frac{P(v)}{2^{n(v)}} \qquad (4.36)$$

**Théorème 4.5.** Soient v un vecteur de parité quelconque et N un entier naturel non nul donc la suite de Collatz de premier terme N et de longueur n(v) admet v comme vecteur de parité si et seulement si $g_v(N) \in \mathbb{N}$ :

$$g_v(N) \in \mathbb{N} \Leftrightarrow V_{n(v)}(N) = v \qquad (4.37)$$

Cette dernière équivalente à :

$$g_v(N) \in \mathbb{N} \Leftrightarrow N \in \mathbb{H}(v)$$

Ou encore :

$$g_v(N) \notin \mathbb{N} \Leftrightarrow N \notin \mathbb{H}(v)$$

Dans ce cas on a:

$$T^{n(v)}(N) = g_v(N) \qquad (4.38)$$

**Démonstration.** Soit N un entier naturel tel que $g_v(N) \in \mathbb{N}$ et M un entier naturel non nul appartient a $\mathbb{H}(v)$ donc la suite de Collatz de premier terme M et de longueur n(v) admet v comme vecteur de parité et dans ce cas on peut écrire :

$$T^{n(v)}(M) = \frac{3^{m(v)}}{2^{n(v)}} M + \frac{P(v)}{2^{n(v)}}$$

Comme on a :

$$g_v(N) = \frac{3^{m(v)}}{2^{n(v)}} N + \frac{P(v)}{2^{n(v)}}$$

Donc on peut écrire :

$$T^{n(v)}(M) - g_v(N) = \frac{3^{m(v)}}{2^{n(v)}} (M - N)$$

Or $T^{n(v)}(M) - g_v(N) \in \mathbb{N}, M - N \in \mathbb{N}$ et $3^{m(v)}$ et $2^{n(v)}$ sont premiers entre eux alors l'équation implique que $(M - N)$ est un multiple de $2^{n(v)}$ c'est à dire que :

$$M - N = 2^{n(v)} j \text{ tel que } j \in \mathbb{Z}$$

Ceci équivalente à :

$$N = M + 2^{n(v)} j$$

Cette dernière égalité implique que $N \in \mathbb{H}(v)$ et par suite la suite de Collatz de premier terme N et de longueur n(v) admet v comme vecteur de parité.

Inversement si N ∈ ℍ(v) donc il est évident que $g_v(N) \in \mathbb{N}$ puisque $T^{n(v)}(N) = g_v(N)$.

**Remarque 4.2.** Soit v un vecteur de parité quelconque donc pour déterminer un entier naturel non nul quelconque appartenant à ℍ(v), il faut trouver un entier naturel non nul tel que son image par la fonction $g_v$ soit un entier naturel.

**Exemple 4.6.** On considère le vecteur de parité v=(1,1,0,1,0) caractérisé par les nombres suivants :

$$\begin{cases} n(v) = 5 \\ m(v) = 3 \\ P(v) = 23 \end{cases}$$

Pour N=5

$$g_v(5) = \frac{3^3}{2^5}x5 + \frac{23}{2^5} = 4.9375 \notin \mathbb{N}$$

La suite de Collatz de premier terme 5 et de longueur 5 n'admet pas v comme vecteur de parité et $5 \notin H(v)$ alors que :

$$g_v(11) = \frac{3^3}{2^5}x11 + \frac{23}{2^5} = 10 \in \mathbb{N}$$

La suite de Collatz de premier terme 11 et de longueur 5 admet v comme vecteur de parité et par suite $11 \in H(v)$.

**Définition 4.7.** Soit v un vecteur de parité de longueur finie, on définit le point caractéristique de v noté X(v) comme le produit de P(v) et a(v) autrement dit :

$$X(v) = P(v)\,a(v) \qquad (4.39)$$

On définit le point caractéristique de v note Y(v) comme le produit de P(v) et b(v) autrement dit :

$$Y(v) = P(v)\,b(v) \qquad (4.40)$$

**Théorème 4.6.** Soit v un vecteur de parité de longueur finie alors le point caractéristique de v noté X(v) est un élément de ℍ(v) c'est à dire que la suite de Collatz de premier terme X(v) et de longueur n(v) admet v comme vecteur de parité ceci se traduit par :

$$X(v) \in \mathbb{H}(v) \qquad (4.41)$$

Donc on peut écrire :

$$T^{n(v)}(X(v)) = \frac{3^{m(v)}}{2^{n(v)}}X(v) + \frac{P(v)}{2^{n(v)}} = Y(v) \qquad (4.42)$$

**Démonstration.** Pour que X(v) soit un élément de ℍ(v) il faut que $g_v(X(v)) \in \mathbb{N}$.

On sait que :

$$g_v(X(v)) = \frac{3^{m(v)}}{2^{n(v)}}X(v) + \frac{P(v)}{2^{n(v)}}$$

$$= \frac{3^{m(v)}}{2^{n(v)}} P(v)\, a(v) + \frac{P(v)}{2^{n(v)}}$$

$$= \frac{P(v)}{2^{n(v)}}\bigl(3^{m(v)}a(v) + 1\bigr)$$

On sait aussi que :

$$3^{m(v)}a(v) + 1 = 2^{n(v)}b(v)$$

On déduit alors que :

$$g_v\bigl(X(v)\bigr) = P(v)b(v) \in \mathbb{N}$$

Donc on peut conclure que :

$$X(v) \in \mathbb{H}(v)$$

Dans ce cas on peut écrire :

$$T^{n(v)}\bigl(X(v)\bigr) = \frac{3^{m(v)}}{2^{n(v)}} P(v)\, a(v) + \frac{P(v)}{2^{n(v)}}$$

$$= P(v)b(v)$$

$$= Y(v)$$

**Exemple 4.7.** On reprend le vecteur de l'exemple précédent défini comme suit :

| 1 | 0 | 1 | 1 | 0 | 1 | 0 | 1 | 1 | 1 |
|---|---|---|---|---|---|---|---|---|---|

Le vecteur de parité considéré est caractérisé par les nombres suivants :

$$\begin{cases} n(v) = 10 \\ m(v) = 7 \\ P(v) = 5645 \end{cases}$$

On sait que :

$$\begin{cases} a(v) = 221 \\ b(v) = 472 \end{cases}$$

Comme P(v)=5645 alors le nombre recherché est le suivant:

$$X(v) = P(v)a(v) = 5645 \times 221 = 1247545$$

On vérifie que la suite de Collatz de premier terme X(v) et de longueur n(v) admet v comme vecteur de parité comme suit:

| 1247545 | 1871318 | 935659 | 1403489 | 2105234 | 1052617 | 1578926 | 789463 |
|---------|---------|--------|---------|---------|---------|---------|--------|
| 1184195 | 1776293 |        |         |         |         |         |        |

**Théorème 4.7.** Soient $(m, n) \in \mathbb{N}\times\mathbb{N}^*$, $v_1$ et $v_2$ deux vecteurs de parité appartenant à l'ensemble $\Omega(n, m)$. On désigne par $x_1 \in \mathbb{H}(v_1)$ et $x_2 \in \mathbb{H}(v_2)$ alors la relation entre $x_1$ et $x_2$ se traduit par l'équation ci-dessous :

$$P(v_1)x_2 - P(v_2)x_1 = 2^n j \text{ tel que } j \in \mathbb{Z} \tag{4.43}$$

**Démonstration.** Comme $(v_1, v_2) \in \Omega(n,m) \times \Omega(n,m)$ celui-ci signifie que :

$$m(v_2) = m(v_1) = m \text{ et } n(v_2) = n(v_1) = n$$

De plus on peut écrire :

$$\begin{cases} \dfrac{3^m x_1 + P(v_1)}{2^n} = y_1 \\ \dfrac{3^m x_2 + P(v_2)}{2^n} = y_2 \end{cases}$$

En multipliant la première équation par $P(v_2)$ et la deuxième par $P(v_1)$, on obtient :

$$\begin{cases} \dfrac{3^m x_1 P(v_2) + P(v_2)P(v_1)}{2^n} = P(v_2)y_1 \\ \dfrac{3^m x_2 P(v_1) + P(v_1)P(v_2)}{2^n} = P(v_1)y_2 \end{cases}$$

L'expression de la différence entre les deux équations s'écrit :

$$\dfrac{3^m}{2^n}\big(x_2 P(v_1) - x_1 P(v_2)\big) = P(v_1)y_2 - P(v_2)y_1$$

On sait que $\big(x_2 P(v_1) - x_1 P(v_2)\big) \in \mathbb{N}$ et $(P(v_1)y_2 - P(v_2)y_1) \in \mathbb{N}$ de plus $3^m$ et $2^n$ sont premiers entre eux donc l'équation ci-dessus signifie que $\big(x_2 P(v_1) - x_1 P(v_2)\big)$ est un multiple de $2^n$.

**Théorème 4.8.** On considère deux vecteurs de parité quelconque $v_1$ et $v_2$ de même longueur n $(n(v_1) = n(v_2)=n)$ de plus on suppose que : $m(v_2) \geq m(v_1)$. Soient $x_1 \in \mathbb{H}(v_1)$, $x_2 \in \mathbb{H}(v_2)$ alors l'expression entre $x_2$ et $x_1$ s'écrit comme suit :

$$3^{m(v_2)-m(v_1)} x_2 P(v_1) - x_1 P(v_2) = 2^n j \text{ tel que } j \in \mathbb{Z} \qquad (4.44)$$

**Démonstration.** On a : $m(v_2) \geq m(v_1)$ et $n(v_2) = n(v_1) = n$

Donc on peut écrire :

$$\begin{cases} \dfrac{3^{m(v_1)} x_1 + P(v_1)}{2^n} = y_1 \\ \dfrac{3^{m(v_2)} x_2 + P(v_2)}{2^n} = y_2 \end{cases}$$

En multipliant la première équation par $P(v_2)$ et la deuxième par $P(v_1)$, on obtient :

$$\begin{cases} \dfrac{3^{m(v_1)} x_1 P(v_2) + P(v_2)P(v_1)}{2^n} = P(v_2)y_1 \\ \dfrac{3^{m(v_2)} x_2 P(v_1) + P(v_1)P(v_2)}{2^n} = P(v_1)y_2 \end{cases}$$

L'expression de la différence entre ces deux équations s'écrit :

$$\dfrac{3^{m(v_1)}}{2^n}\Big(3^{m(v_2)-m(v_1)} x_2 P(v_1) - x_1 P(v_2)\Big) = P(v_1)y_2 - P(v_2)y_1$$

On sait que $\left(3^{m(v_2)-m(v_1)}x_2 P(v_1) - x_1 P(v_2)\right) \in \mathbb{N}$ et $(P(v_1)y_2 - P(v_2)y_1) \in \mathbb{N}$ de plus $3^{m(v_1)}$ et $2^n$ sont premiers entre eux donc l'équation ci-dessus signifie que :

$\left(3^{m(v_2)-m(v_1)}x_2 P(v_1) - x_1 P(v_2)\right)$ est un multiple de $2^n$.

**Définition 4.8.** Soient n un entier naturel non nul, m un entier naturel quelconque et v un vecteur de parité finie tel que $n(v) = n$ et $m(v) = m$.

$$P(v) = 3^{m-1}x2^{j_1-1} + 3^{m-2}x2^{j_2-1} + \cdots + 3x2^{j_{m-1}-1} + 2^{j_m-1}$$

Pour tout entier naturel non nul k allant de 1 à n, on définit les deux nombres $\theta_k(v)$ et $t_k(v)$ vérifiant l'équation suivante :

$$3^k \theta_k(v) + 1 = 2^{n-j_k+1} t_k(v) \tag{4.45}$$

$(\theta_k(v), t_k(v)) \in \mathbb{N}^* x \mathbb{N}^*$

En multipliant l'équation précédente par $2^{j_k-1}$, on obtient la deuxième équation associée au nombre P(v) définie comme suit :

$$3^k z_k(v) + 2^{j_k-1} = 2^n t_k(v) \tag{4.46}$$

Avec :

$$z_k(v) = 2^{j_k-1} \theta_k(v) \tag{4.47}$$

En fin on définit les deux nombres $\mathbf{X}^*(\mathbf{v})$ et $\mathbf{Y}^*(\mathbf{v})$ comme suit :

$$\begin{cases} \mathbf{X}^*(\mathbf{v}) = \sum_{k=1}^{m} z_k(v) = \sum_{i=1}^{m} 2^{j_k-1} \theta_k(v) \\ \mathbf{Y}^*(\mathbf{v}) = \sum_{k=1}^{m} 3^{m-i} t_k(v) \end{cases} \tag{4.48}$$

**Théorème 4.9.** On considère un vecteur de parité **de longueur finie** n défini comme suit :

$$v = (e_1, e_2, \ldots, e_n)$$

Tel que $\forall\ 1 \leq i \leq n$, $e_i \in \{0,1\}$

On suppose que ce vecteur contient m termes qui sont égaux à 1. On suppose aussi que ces termes sont distribués de la façon suivante dans le vecteur considéré:

$$v = \left(0,0,\ldots,0, e_{j_1}, 0,0\ldots,0, e_{j_2}, 0,0,\ldots,0, e_{j_i}, 0\ldots,0, e_{j_m}, 0,\ldots,0\right) \mid \forall\ 1 \leq i \leq m, e_{j_i} = 1$$

L'expression du nombre caractéristique P(v) du vecteur v en fonction de m et de $j_1, j_2, \ldots, j_m$ s'écrit comme suit :

$$\mathbf{P(v) = 3^{m-1}x2^{j_1-1} + 3^{m-2}x2^{j_2-1} + 3^{m-3}2^{j_3-1} + \cdots + 3x2^{j_{m-1}-1} + 2^{j_m-1}}$$

**Donc le nombre $\mathbf{X}^*(\mathbf{v})$ qui s'écrit sous la forme ci-dessous est un élément $\mathbf{X}^* \in \mathbb{H}(\mathbf{v})$ qui peut s'écrire sous la forme suivante :**

$$X^*(v) = \sum_{k=1}^{m} 2^{j_k - 1} \theta_k(v) = 2^{j_1 - 1} \theta_1(v) + 2^{j_2 - 1} \theta_2(v) + \cdots + 2^{j_m - 1} \theta_m(v) \quad (4.49)$$

Tel que $\theta_1(v), \theta_2(v), \ldots, \theta_m(v)$ sont des entiers naturels impairs autrement $\forall\ 1 \leq i \leq m$, on a $\theta_i \in (2\mathbb{N} + 1)$, de plus on a :

$$T^{n(v)}(X^*(v)) = Y^*(v) \quad (4.50)$$

**Démonstration.** L'expression de $3^m X^*(v) + P(v)$ est la suivante :

$$3^m X^*(v) + P(v) = 3^m \sum_{i=1}^{m} z_k + P(v) =$$

$$3^{m-1} x (3z_1 + 2^{j_1 - 1}) + 3^{m-2} x (3^2 z_2 + 2^{j_2 - 1}) + 3^{m-3}(3^3 z_3 + 2^{j_3 - 1}) + \cdots + 3x(2^{j_{m-1} - 1} + 3^{m-1} z_{m-1}) + 3^m z_m + 2^{j_m - 1}$$

Ou encore :

$$3^m X^*(v) + P(v) = \sum_{k=1}^{m} 3^{m-k}(3^k z_k(v) + 2^{j_k - 1})$$

$\forall\ k \in [\![1, m]\!]$, on sait que :

$$3^k z_k(v) + 2^{j_k - 1} = 2^n \delta_k(v)$$

Dans ce cas chaque terme qui s'écrit sous la forme $3^k z_k(v) + 2^{j_k - 1}$ est divisible par $2^n$ et dans ce cas on peut écrire :

$$3^m X^*(v) + P(v) = 2^n \sum_{k=1}^{m} 3^{m-k} t_k(v)$$

Ceci équivalent à :

$$\frac{3^m}{2^n} X^*(v) + \frac{P(v)}{2^n} = \sum_{k=1}^{m} 3^{m-k} t_k(v) = Y^*(v)$$

Ce qui implique que :

$$g_v(X^*(v)) = \sum_{k=1}^{m} 3^{m-k} t_k(v) \in \mathbb{N}$$

D'après le théorème précédent $X^*(v)$ est un élément de $\mathbb{H}(v)$

L'expression de $X^*(v)$ s'écrit comme suit:

$$X^*(v) = \sum_{k=1}^{m} 2^{j_k - 1} \theta_k(v) = 2^{j_1 - 1} \theta_1(v) + 2^{j_2 - 1} \theta_2(v) + \cdots + 2^{j_m - 1} \theta_m(v)$$

$$\forall\ k \in [\![1, m]\!] \text{ on a } \theta_k(v) \in 2\mathbb{N} + 1$$

**Exemple 4.8.** On considère le vecteur de parité ci-dessous de longueur $n(v) = 10$:

| 1 | 0 | 1 | 1 | 0 | 1 | 0 | 1 | 1 | 1 |
|---|---|---|---|---|---|---|---|---|---|

Les deux nombres caractéristiques de vecteur v sont les suivant:
$$\begin{cases} n(v) = 10 \\ m(v) = 7 \end{cases}$$

Le nombre caractéristiques $P(v)$ est donné par:
$$P(v) = 3^6 + 3^5 \, x2^2 + 3^4 \, x2^3 + 3^3 x \, 2^5 + 3^2 x 2^7 + 3 \, x2^8 + 2^9$$

On cherche à déterminé un entier naturel non nul qu'on va le noter $X^*(v)$ qui peut s'écrire sous la forme suivante:

$$X^*(v) = \sum_{k=1}^{m} 2^{j_k - 1} \theta_k = \theta_1 + 2^2 \theta_2 + 2^3 \theta_3 + 2^5 \theta_4 + 2^7 \theta_5 + 2^8 \theta_6 + 2^9 \theta_7$$

Noter que :
$$\frac{3^{m(v)} X^*(v) + P(v)}{2^{n(v)}} \in \mathbb{N}$$

On suppose que:
$$X^*(v) = \sum_{i=1}^{m(v)} z_k(v)$$

Pour tout entier naturel non nul j allant de 1 à 7, on doit déterminer les couples $(z_j, t_j)$ tel que:
$$2^{j_i - 1} + 3^j z_j(v) = 2^{10} t_j(v)$$

Donc on a un système à 7 équations comme suit :
$$\begin{cases} 3z_1 + 1 = 2^{10} t_1 \\ 3^2 z_2 + 2^2 = 2^{10} t_2 \\ 3^3 z_3 + 2^3 = 2^{10} t_3 \\ 3^4 z_4 + 2^5 = 2^{10} t_4 \\ 3^5 z_5 + 2^7 = 2^{10} t_5 \\ 3^6 a_6 + 2^8 = 2^{10} t_6 \\ 3^7 z_7 + 2^9 = 2^{10} t_7 \end{cases}$$

Pour le cas du couple $(z_1, t_1)$, on ecrit:
$$1 + 3z_1 = 2^{10} t_1$$

On remarque que:
$$1 + 3x341 = 2^{10}$$

Il en result que:
$$\begin{cases} z_1 = 341 \\ t_1 = 1 \end{cases}$$

Pour le cas du couple $(z_2, t_2)$, on écrit:

$$2^2 + 3^2 z_2 = 2^{10} t_2$$

L'équation recherchée est la suivante:

$$2^2 + 3^2 \text{x} 796 = 2^{10} \text{x} 7$$

On déduit alors que:

$$\begin{cases} z_2 = 796 \\ t_2 = 7 \end{cases}$$

Pour le cas du couple $(z_3, t_3)$, on écrit:

$$2^3 + 3^3 z_3 = 2^{10} t_3$$

D'après le calcul

$$2^3 + 3^3 \text{x} 872 = 2^{10} \text{x} 23$$

Cette dernière équation nous permet d'écrire :

$$\begin{cases} z_3 = 872 \\ t_3 = 23 \end{cases}$$

Pour le cas du couple $(z_4, t_4)$, on écrit:

$$2^5 + 3^4 z_4 = 2^{10} t_4$$

L'équation qui correspond :

$$2^5 + 3^4 \text{x} 480 = 2^{10} \text{x} 38$$

Ceci signifie que :

$$\begin{cases} z_4 = 480 \\ t_4 = 38 \end{cases}$$

Pour le cas du couple $(z_5, t_5)$, on écrit:

$$2^7 + 3^5 z_5 = 2^{10} t_5$$

Les deux nombres recherchés sont les suivants :

$$2^7 + 3^5 640 = 2^{10} \text{x} 152$$

Pour le cas du couple $(z_6, t_6)$, on écrit:

$$2^8 + 3^6 z_6 = 2^{10} t_6$$

$$2^8 + 3^6 \text{x} 768 = 2^{10} \text{x} 547$$

Pour le cas du couple $(z_7, t_7)$, on écrit:

$$3^7 z_7 + 2^9 = 2^{10} t_7$$

$$3^7 2^9 + 2^9 = 2^{10} \text{x} 1094$$

$$z_1 = 341, z_2 = 796, z_3 = 872, z_4 = 480, z_5 = 640, z_6 = 768, z_7 = 512$$

$$X^*(v) = \sum_{i=1}^{7} z_i = 341 + 796 + 872 + 480 + 640 + 768 + 512 = 4409$$

$$Y^*(v) = \sum_{i=1}^{7} 3^{m(v)-i} t_i$$

$$= 3^6 x1 + 3^5 x7 + 3^4 x23 + 3^3 x38 + 3^2 x152 + 3x547 + 1094 = 9422$$

On sait que:

$$\theta_i = \frac{z_i}{2^{j_i}}$$

On déduit les valeurs des $\theta_i$ $\forall$ i allant de 1 à 7 comme suit:

$$\theta_1 = 341, \theta_2 = 199, \theta_3 = 109, \theta_4 = 15, \theta_5 = 5, \theta_6 = 3, \theta_7 = 1$$

L'expression de x en forme polynomiale est la suivante:

$$X^*(v) = 341 + 2^2 x199 + 2^3 x109 + 2^5 x15 + 2^7 x5 + 2^8 x3 + 2^9$$

La suite de Collatz de premier terme $X^*(v)$ et de longueur 10 est la suivante:

| 4409 | 6614 | 3307 | 4961 | 7442 | 3721 | 5582 | 2791 | 4187 | 6281 |

On déduit aussi:

$$X^*(v) = N_4(v)$$

En effet:

$$4409 = 313 + 2^{10} x4$$

| $N_0(v)$ | $N_1(v)$ | $N_2(v)$ | $N_3(v)$ | $N_4(v)$ | | $N_j(v)$ | |
|---|---|---|---|---|---|---|---|
| 313 | 1337 | 2361 | 3385 | 4409 | | $313 + 2^{10}j$ | |

**Corollaire 4.3.** Soit v un vecteur de parité de longueur infini, donc la relation entre les deux points particuliers $X(v)$ et $X^*(v)$ est comme suit :

$$X(v) = X^*(v) + 2^{n(v)} J(v) \qquad (4.51)$$

Avec $J(v)$ un entier relatif.

## 5. Caractérisation des vecteurs de parité de longueur infinie $n(v) = \infty$

Soit V un vecteur de parité de longueur infinie qu'on peut le représenter comme suit :

$$V = (e_1, e_2, \dots, e_j, \dots) \qquad (5.1)$$

Tel que $\forall j \in \mathbb{N}^*$ $e_j \in \{0,1\}$

Ce vecteur est caractérisé par une longueur infinie donc on peut écrire :

$$n(v) = \infty \qquad (5.2)$$

Il est inutile de considérer de vecteur de parité de longueur infinie mais son nombre caractéristique m(v) est fini donc dans le cas de l'étude des vecteurs de parité de longueurs infinies on suppose toujours qu'ils contiennent un nombre infini des termes qui sont égaux à 1 autrement dit :

$$m(v) = \infty \tag{5.3}$$

Rappeler qu'une partie de V qui contient les j premiers éléments de V est notée comme suit :

$$R_j(V) = (e_1, e_2, \ldots, e_j)$$

Rappeler aussi que pour un V un vecteur de parité de longueur infinie V on peut écrire :

$$\lim_{j \to +\infty} R_j(V) = V$$

**Notations 5.1.** Soit j un entier naturel non nul ($j \in \mathbb{N}^*$) alors les principaux nombres caractéristiques d'ordre j du vecteur V sont définies comme suit :

$$n\left(R_j(V)\right) = n_j(V)$$

$$m\left(R_j(V)\right) = m_j(V)$$

$$P\left(R_j(V)\right) = P_j(V)$$

$$c\left(R_j(V)\right) = c_j(V)$$

Les points particuliers d'ordre j du vecteur V sont comme suit :

$$X\left(R_j(V)\right) = X_j(V)$$

$$Y\left(R_j(V)\right) = Y_j(V)$$

$$X^*\left(R_j(V)\right) = X_j^*(V)$$

$$Y^*\left(R_j(V)\right) = Y_j^*(V)$$

**Définition 5.1.** L'équation caractéristique d'ordre j du vecteur binaire V est la suivante :

$$3^{m_j(V)} a_j(V) + 1 = 2^{n_j(V)} b_j(V) \tag{5.4}$$

Avec :

$$a_j(V) = a\left(R_j(V)\right), \quad b_j(V) = b\left(R_j(V)\right)$$

**Notation 5.2.** Si on considère une suite de Collatz de premier terme N et de vecteur de parité $R_j(V)$ donc l'expression de $T^{n_j(V)}(N)$ en fonction de N et les nombres caractéristiques principaux d'ordre j de V s'écrit comme suit :

$$T^{n_j(V)}(N) = \frac{3^{m_j(V)}}{2^{n_j(V)}}N + \frac{P_j(V)}{2^{n_j(V)}} \quad (5.5)$$

**Définition 5.2.** On définit une deuxième équation caractéristique d'ordre j pour le vecteur binaire infini V en multipliant l'équation par $P_j(V)$, on obtient :

$$3^{m_j(V)}a_j(V)P_j(V) + P_j(V) = 2^{n_j(V)}b_j(V)P_j(V)$$

Cette dernière équation peut être écrite sous différentes formes par exemple si on devise cette dernière par $2^{n_j(V)}$ on obtient :

$$\frac{3^{m_j(V)}}{2^{n_j(V)}}a_j(V)P_j(V) + \frac{P_j(V)}{2^{n_j(V)}} = b_j(V)P_j(V)$$

Comme $a_j(V)P_j(V) = X_j(V)$ et $b_j(V)P_j(V) = Y_j(V)$ donc on peut écrire :

$$\frac{3^{m_j(V)}}{2^{n_j(V)}}X_j(V) + \frac{P_j(V)}{2^{n_j(V)}} = Y_j(V) \quad (5.6)$$

**Corollaire 5.1.** Soit V un vecteur de parité de longueur infinie donc pour tout entier naturel non nul j, La suite de Collatz de premier terme $X_j(V) = a_j(V)P_j(V)$ de longueur $n_j(V)$ admet $R_j(V)$ comme vecteur de parité donc on peut écrire :

$$T^{n_j(V)}(X_j(V)) = \frac{3^{m_j(V)}}{2^{n_j(V)}}X_j(V) + \frac{P_j(V)}{2^{n_j(V)}}$$

**Corollaire 5.2.** Soit V un vecteur de parité quelconque de longueur infinie alors pour tout entier naturel non nul j, il existe un entier naturel non nul qu'on le note $K_j(V)$ tel que la relation entre $X_j(V)$ et $N_{0,j}(V)$ peut s'écrire sous la forme suivante :

$$X_j(V) = N_{0,j}(V) + 2^{n_j(V)}K_j(v) \quad \text{avec} \quad K_j(v) \in \mathbb{N}^* \quad (5.7)$$

Noter que $K_j(V)$ figure parmi les nombres caractéristiques d'ordre j du vecteur V. de plus il existe un nombre caractéristique de v noté $K_j^*(V) \in \mathbb{N}^*$ tel que :

$$X_j^*(V) = N_{0,j}(V) + 2^{n_j(V)}K_j^*(V) \quad (5.8)$$

**Démonstration.** On sait bien que $X_j(V) \in \mathbb{H}_j(V)$ et comme $N_{0,j}(V) \in \mathbb{H}_j(V)$ donc il est évident que la différence entre ces deux nombres est un multiple de $2^{n_j(V)}$. Noter que $X_j(V) > N_{0,j}(V)$ donc on peut conclure qu'il existe un entier naturel $K_j(V)$ caractéristique de V tel que :

$$X_j(V) = N_{0,j}(V) + 2^{n_j(V)}K_j(V)$$

De même pour $X_j^*(V)$ qui appartient aussi au $\mathbb{H}_j(V)$

**Définition 5.3.** On définit le rapport caractéristique d'ordre j d'un vecteur binaire de longueur infinie V comme le rapport de son point caractéristique d'ordre j noté $X_j(V)$ par $2^{n_j(V)}$ qui s'écrit comme suit :

$$q_j(V) = \frac{X_j(V)}{2^{n_j(V)}} \tag{5.9}$$

De la même manière, on définit le rapport $q_j^*(V)$ caractéristique comme suit :

$$q_j^*(V) = \frac{X_j^*(V)}{2^{n_j(V)}}$$

**Définition 5.4.** On considère un vecteur binaire quelconque V de longueur infinie. On sait que pour tout entier naturel non nul j, le vecteur de parité $R_j(V)$ est caractérisée par un nombre rationnel noté $r_0\big(R_j(V)\big)$ qui représente le rapport de $N_0(R_j(V))$ par $2^{n(R_j(V))}$. Ce même nombre caractéristique représente un nombre caractéristique d'ordre j pour le vecteur binaire infini V. Ce nombre s'écrit comme suit :

$$r_{0,j}(V) = \frac{N_{0,j}(V)}{2^{n_j(V)}} \tag{5.10}$$

**Corollaire 5.3.** Soit V un vecteur de parité de longueur infinie donc pour tout entier naturel non nul j, le vecteur infinie V peut être caractérisé par l'équation caractéristique d'ordre j suivante :

$$q_j(V) = r_{0,j}(V) + K_j(V) \tag{5.11}$$

Noter que :

$$\big(q_j(V), r_{0,j}(V)\big) \in \mathbb{R}^* x \mathbb{R}^*, K_j(V) \in \mathbb{N}$$

**Notation 5.3.** La limite de $N_{0,j}(V)$ lorsque j tend vers l'infini est notée comme suit:

$$\lim_{j \to +\infty} N_{0,j}(V) = N_{0,\infty}(V)$$

**Lemme 5.1.** Soit V un vecteur binaire de longueur infinie et pour un entier naturel non nul j quelconque on considère les deux vecteurs binaires $R_j(V)$ et $R_{j+1}(V)$ donc on peut distinguer deux cas possibles pour la relation entre $r_{0,j}(V)$ et $r_{0,j+1}(V)$ :

**Premier cas**: $N_{0,j}(V) = N_{0,j+1}(V)$

Dans ce cas les deux suites de Collatz notées $S_j(V)$ et $S_{j+1}(V)$ possèdent le même premier terme et on peut écrire alors :

$$\frac{N_{0,j+1}(V)}{2^{n_{j+1}(V)}} = \frac{1}{2} \frac{N_{0,j}(V)}{2^{n_j(V)}} \tag{5.12}$$

Ce qui implique que :

$$r_{0,j+1}(V) = \frac{1}{2} r_{0,j}(V)$$

**Deuxième cas: $N_{0,j}(V) \neq N_{0,j+1}(V)$**

Les deux suites $S_j(V)$ et $S_{j+1}(V)$ ne possèdent pas les mêmes premiers termes et la relation entre leurs premiers termes est comme suit :

$$N_{0,j+1}(V) = N_{0,j}(V) + 2^{n_j(V)}$$

Ceci nous permet d'écrire :

$$\frac{N_{0,j+1}(V)}{2^{n_{j+1}(V)}} = \frac{1}{2} \frac{N_{0,j}(V)}{2^{n_j(V)}} + \frac{1}{2}$$

Il en résulte que :

$$r_{0,j+1}(V) = \frac{1}{2} r_{0,j}(V) + \frac{1}{2} \tag{5.13}$$

**Définition 5.5.** Les nombres caractéristiques absolus d'un vecteur de parité de longueur infinie sont les limites des nombres caractéristiques d'ordre finie autrement si on désigne par $x_j(V)$ un nombre caractéristique d'ordre j de vecteur V alors le nombre caractéristique absolus de V noté $x(V)$ est défini comme suit :

$$x(V) = \lim_{j \to +\infty} x_j(V)$$

Le nombre caractéristique $P(V)$ est défini comme suit:

$$P(V) = \lim_{j \to +\infty} P_j(V) \tag{5.14}$$

Le nombre de Collatz caractéristique de vecteur de parité V est le suivant :

$$c(V) = \lim_{j \to +\infty} c_j(V) \tag{5.15}$$

**Lemme 5.2.** Soit V un vecteur binaire de longueur infinie. On sait La suite constituée par les premiers termes des suites $S_{0,1}(V), S_{0,2}(V), \ldots, S_{0,j}(V), \ldots$ est notée comme suit :

$$h_0(V) = \left(N_{0,1}(V), N_{0,2}(V), N_{0,3}(V), \ldots, N_{0,j}(V), \ldots\right)$$

Si on a :

$$\lim_{j \to +\infty} N_{0,j}(V) = +\infty$$

Alors elle existe une suite de longueur infinie extraite de la suite $h_0(V)$ qui contient des éléments **deux à deux distincts**. Cette suite extraite est notée comme ci-dessous :

$$h_{ext,0}(V) = \left(N_{0,k_1}(V), N_{0,k_2}(V), N_{0,k_3}(V), \ldots, N_{0,k_i}(V), \ldots\right) \tag{5.16}$$

Tel que :

$$\forall (i,j) \in \mathbb{N}^* \times \mathbb{N}^* \text{ si } i \neq j \text{ alors } N_{0,k_i}(V) \neq N_{0,k_j}(V) \tag{5.17}$$

Dans ce cas on peut écrire :

$$N_{0,k_1}(V) < N_{0,k_2}(V) < N_{0,k_3}(V) < \cdots < N_{0,k_i}(V) < \cdots \tag{5.18}$$

**Définition 5.6.** Soit un **vecteur de parité V de longueur infinie**, on définit les deux nombres caractéristiques $\alpha_j(V)$ et $\beta_j(V)$ qui sont deux entiers naturels non nuls tel que :

$$P_j(V) = 3^{m_j(V)}\alpha_j(V) + \beta_j(V) \mid \beta_j(V) < 3^{m_j(V)} \tag{5.19}$$

De plus, on définit les deux nombres caractéristiques $A_j(V)$ et $B_j(V)$ qui sont deux entiers naturels non nuls tel que :

$$P_j(V) = 2^{n_j(V)}A_j(V) + B_j(V) \mid B_j(V) < 2^{n_j(V)} \tag{5.20}$$

**Corollaire 5.4** Soit V un vecteur binaire de longueur infinie alors pour tout entier naturel non nul j, le nombre caractéristique d'ordre j de vecteur V noté $P_j(V)$ vérifie l'inégalité suivante :

$$3^{m(v_j)} - 2^{m(v_j)} \leq P_j(V) \leq 2^{n(v_j)-m(v_j)}\left(3^{m(v_j)} - 2^{m(v_j)}\right) \tag{5.21}$$

**Corollaire 5.6** Soit V un vecteur binaire de longueur infinie alors on peut écrire :

$$\lim_{j \to +\infty} \frac{P_j(V)}{2^{n_j(V)}3^{m_j(V)}} = 0 \tag{5.22}$$

Noter qu'on peut remplacer $j \to +\infty$ par $n_j(V) \to +\infty$

**Démonstration.** On sait que :

$$P_j(V) \leq 2^{n(v_j)-m(v_j)}(3^{m(v_j)} - 2^{m(v_j)})$$

On fait diviser inégalité de deux cotes par $2^{n_j(V)}3^{m_j(V)}$, on obtient :

$$\frac{P_j(V)}{2^{n_j(V)}3^{m_j(V)}} \leq \frac{1}{2^{m_j(V)}} - \frac{1}{3^{m_j(V)}}$$

Comme

$$\lim_{j \to +\infty} \left(\frac{1}{2^{m_j(V)}} - \frac{1}{3^{m_j(V)}}\right) = 0$$

Donc on peut déduire que :

$$\lim_{j \to +\infty} \left(\frac{P_j(V)}{2^{n_j(V)}3^{m_j(V)}}\right) = 0$$

**Corollaire 5.7.** Soit V un vecteur binaire de longueur infinie. La limite de rapport de son nombre caractéristique $P_j(V)$ d'ordre j par $3^{m_j(V)}$ satisfait la condition suivante :

$$\lim_{j \to +\infty} \frac{P_j(V)}{3^{m_j(V)}} \geq 1 \tag{5.23}$$

**Démonstration.** On sait que :

$$P_j(V) \geq (3^{m(v_j)} - 2^{m(v_j)})$$

Ce qui implique :

$$P_j(V) \geq (3^{m(v_j)})(1 - \frac{2^{m(v_j)}}{3^{m(v_j)}})$$

On divise par $3^{m_j(V)}$, on obtient :

$$\frac{P_j(V)}{3^{m_j(V)}} \geq 1 - \frac{2^{m(v_j)}}{3^{m(v_j)}}$$

Ce qui implique :

$$\lim_{j \to +\infty} \frac{P_j(V)}{3^{m_j(V)}} \geq \lim_{j \to +\infty} (1 - \frac{2^{m(v_j)}}{3^{m(v_j)}})$$

Comme

$$\lim_{j \to +\infty} (1 - \frac{2^{m(v_j)}}{3^{m(v_j)}}) = 1$$

Donc on déduit que :

$$\lim_{j \to +\infty} \frac{P_j(V)}{3^{m_j(V)}} \geq 1$$

**Corollaire 5.8.** Soit V un vecteur binaire de longueur infinie. Pour tout entier naturel non nul j, l'expression de $P_j(V)$ en fonction de $3^{m_j(V)}$ s'écrit comme suit :

$$P_j(V) = 3^{m_j(V)} \alpha_j(V) + \beta_j(V) \text{ avec } \beta_j(V) < 3^{m_j(V)}$$

Alors le rapport de $\alpha_j(V)$ par $2^{n_j(V)}$ admet la limite suivante :

$$\lim_{j \to +\infty} \frac{\alpha_j(V)}{2^{n_j(V)}} = 0 \tag{5.24}$$

C'est à dire que pour des longueurs suffisamment grandes de $\alpha_j(V)$ devient négligeable par rapport à $2^{n_j(V)}$. Autrement **au voisinage de l'infini**, le coefficient $\alpha_j(V)$ est négligeable devant $2^{n_j(V)}$.

$$\alpha_j(V) = \circ \left(2^{n_j(V)}\right) \tag{5.25}$$

**Démonstration.** L'expression de $P_j(V)$ en fonction de $3^{m_j(V)}$ est comme suit :

$$P_j(V) = 3^{m_j(V)} \alpha_j(V) + \beta_j(V)$$

On sait que :

$$\frac{P_j(V)}{2^{n_j(V)} 3^{m_j(V)}} \leq \frac{1}{2^{m_j(V)}} - \frac{1}{3^{m_j(V)}}$$

Ce qui implique :

$$\frac{3^{m_j(V)} \alpha_j(V) + \beta_j(V)}{2^{n_j(V)} 3^{m_j(V)}} \leq \frac{1}{2^{m_j(V)}} - \frac{1}{3^{m_j(V)}}$$

Ou encore :

$$\left(\frac{\alpha_j(V)}{2^{n_j(V)}} + \frac{\beta_j(V)}{2^{n_j(V)}3^{m_j(V)}}\right) \leq \frac{1}{2^{m_j(V)}} - \frac{1}{3^{m_j(V)}}$$

Les vecteurs binaires de longueurs infinies qui font l'objet de l'étude sont caractérisés par des nombres $m_j(V)$ infinis donc on peut écrire :

$$\lim_{j \to +\infty} \left(\frac{1}{2^{m_j(V)}} - \frac{1}{3^{m_j(V)}}\right) = 0$$

On déduit que :

$$\lim_{j \to +\infty} \left(\frac{\alpha_j(V)}{2^{n_j(V)}} + \frac{\beta_j(V)}{2^{n_j(V)}3^{m_j(V)}}\right) = 0$$

De plus on a :

$$\frac{\beta_j(V)}{2^{n_j(V)}3^{m_j(V)}} < \frac{P_j(V)}{2^{n_j(V)}3^{m_j(V)}}$$

On déduit que :

$$\lim_{j \to +\infty} \left(\frac{\beta_j(V)}{2^{n_j(V)}3^{m_j(V)}}\right) = \lim_{j \to +\infty} \left(\frac{P_j(V)}{2^{n_j(V)}3^{m_j(V)}}\right) = 0$$

On peut conclure que :

$$\lim_{j \to +\infty} \left(\frac{\alpha_j(V)}{2^{n_j(V)}}\right) = 0$$

**Corollaire 5.9.** Soit V un vecteur binaire de longueur infinie. Pour tout entier naturel non nul j, l'expression de $P_j(V)$ en fonction de $2^{n_j(V)}$ s'écrit comme suit :

$$P_j(V) = 2^{n_j(V)}A_j(V) + B_j(V) \text{ avec } B_j(V) < 2^{n_j(V)}$$

Alors le rapport de $A_j(V)$ par $3^{m_j(V)}$ admet la limite suivante :

$$\lim_{j \to +\infty} \frac{A_j(V)}{3^{m_j(V)}} = 0 \tag{5.26}$$

Autrement **au voisinage de l'infini**, le coefficient $A_j(V)$ est négligeable devant $3^{m_j(V)}$.

$$A_j(V) = \circ\left(3^{m_j(V)}\right) \tag{5.27}$$

**Démonstration.** On sait que :

$$\frac{P_j(V)}{2^{n_j(V)}3^{m_j(V)}} \leq \frac{1}{2^{m_j(V)}} - \frac{1}{3^{m_j(V)}}$$

Ce qui implique :

$$\frac{2^{n_j(V)}A_j(V) + B_j(V)}{2^{n_j(V)}3^{m_j(V)}} \leq \frac{1}{2^{m_j(V)}} - \frac{1}{3^{m_j(V)}}$$

Ou encore :

$$\left(\frac{A_j(V)}{3^{m_j(V)}} + \frac{B_j(V)}{2^{n_j(V)}3^{m_j(V)}}\right) \leq \frac{1}{2^{m_j(V)}} - \frac{1}{3^{m_j(V)}}$$

Les vecteurs binaires de longueurs infinies qui font l'objet de l'étude sont caractérisés par des nombres $m_j(V)$ infinis donc on peut écrire :

$$\lim_{j \to +\infty} \left(\frac{1}{2^{m_j(V)}} - \frac{1}{3^{m_j(V)}}\right) = 0$$

On déduit que :

$$\lim_{j \to +\infty} \left(\frac{A_j(V)}{3^{m_j(V)}} + \frac{B_j(V)}{2^{n_j(V)}3^{m_j(V)}}\right) = 0$$

De plus on a :

$$\frac{B_j(V)}{2^{n_j(V)}3^{m_j(V)}} < \frac{P_j(V)}{2^{n_j(V)}3^{m_j(V)}}$$

On déduit que :

$$\lim_{j \to +\infty} \left(\frac{B_j(V)}{2^{n_j(V)}3^{m_j(V)}}\right) = \lim_{j \to +\infty} \left(\frac{P_j(V)}{2^{n_j(V)}3^{m_j(V)}}\right) = 0$$

On peut conclure que :

$$\lim_{j \to +\infty} \left(\frac{A_j(V)}{3^{m_j(V)}}\right) = 0$$

**Corollaire 5.10.** Soit V un vecteur binaire infini donc les deux suites $(X_j(V)/2^{n_j(V)})_{j \geq 1}$ et $(Y_j(V)/3^{m_j(V)})_{j \geq 1}$ sont équivalente au voisinage de l'infini autrement dit :

$$\frac{X_j(V)}{2^{n_j(V)}} \sim \frac{Y_j(V)}{3^{m_j}}$$

Cette dernier équivaut à :

$$\frac{a_j(V)}{2^{n_j(V)}} \sim \frac{b_j(V)}{3^{m_j}}$$

## 6. Conditions d'existence et de non existence des suites divergente

**Lemme 6.1.** Soit V un vecteur binaire de longueur infinie. Concernant, le comportement asymptotique du $r_{0,j}(V)$ au voisinage de l'infini, on peut distinguer deux cas différents :

**Premier cas :**

Le coefficient $r_{0,j}(V)$ tend vers 0 lorsque j tend vers l'infini c'est à dire que V est caractérisé par un rapport $r_{0,\infty}(V)$ pratiquement nul.

$$r_{0,\infty}(V) = \lim_{j \to +\infty} \frac{N_{0,j}(V)}{2^{n_j(V)}} = 0 \tag{6.1}$$

Dans ce cas, on montre que $N_{0,j}(V)$ admet nécessairement une limite finie à l'infini.

**Deuxième cas :**

Le coefficient $r_{0,j}(V)$ ne tend par vers 0 lorsque j tend vers l'infini c'est à dire que V est caractérisé par un rapport absolu $r_{0,\infty}(V)$ non nul qui s'écrit :

$$\lim_{j\to+\infty} \frac{N_{0,j}(V)}{2^{n_j(V)}} \neq 0 \qquad (6.2)$$

Dans ce cas il est évident que $N_{0,j}(V)$ ne peut pas admettre une limite finie elle tend nécessairement vers l'infini.

**Théorème 6.1.** Soit V un vecteur de parité de longueur infinie et on considère les deux suites caractéristiques $(N_{0,j}(V))_{j\in\mathbb{N}^*}$ et $(r_{0,j}(V))_{j\in\mathbb{N}^*}$ relatives au vecteur infini V donc la suite $(r_{0,j}(V))_{j\in\mathbb{N}^*}$ admet une limite infini à l'infini **si et seulement si** la suite $(r_{0,j}(V))_{j\in\mathbb{N}^*}$ n admet pas 0 comme limite lorsque j tend vers l'infini, autrement dit :

$$\lim_{j\to+\infty} \frac{N_{0,j}(V)}{2^{n_j(V)}} \neq 0 \Leftrightarrow \lim_{j\to+\infty} N_{0,j}(V) = +\infty \qquad (6.3)$$

**Démonstration.** On suppose que les deux conditions suivantes sont satisfaites en même temps :

$$\lim_{j\to+\infty} \frac{N_{0,j}(V)}{2^{n_j(V)}} = 0 \text{ et } \lim_{j\to+\infty} N_{0,j}(V) = +\infty$$

Si on a :

$$\lim_{j\to+\infty} \frac{N_{0,j}(V)}{2^{n_j(V)}} = 0$$

Donc on peut écrire :

$$\forall \varepsilon > 0, \exists j_0 \text{ tel que } \forall j > j_0 \text{ on a } \frac{N_{0,j}(V)}{2^{n_j(V)}} < \varepsilon$$

D'autre part comme si :

$$\lim_{j\to+\infty} N_{0,j}(V) = +\infty$$

Donc $\forall j \in \mathbb{N}^*$ il existe un entier naturel non nul k ($k \geq 1$) tel que :

$$N_{0,j}(V) \neq N_{0,j+k}(V)$$

En particulier, il existe un entier naturel non nul k tel que :

$$N_{0,j_0}(V) \neq N_{0,j_0+k}(V)$$

Une telle condition, on peut l'exprimer comme suit :

$$N_{0,j_0}(V) = N_{0,j_0+1}(V) = \cdots = N_{0,j_0+k-1}(V) \neq N_{0,j_0+k}(V)$$

On prend un entier naturel j tel que $j > j_0$ on sait que :

$$\frac{N_{0,j_0}(V)}{2^{n_j(V)}} < \varepsilon$$

Comme $N_{0,j_0}(V) = N_{0,j_0+1}(V) = \cdots = N_{0,j_0+k-1}(V)$ alors $\forall$ i tel que $0 \leq 1 \leq k-1$

$$\frac{N_{0,j_0+i}(V)}{2^{n_{j+i}(V)}} < \varepsilon$$

Comme $N_{0,j_0+k-1}(V) \neq N_{0,j_0+k}(V)$ donc nécessairement :

$$N_{0,j_0+k}(V) = N_{0,j_0+k-1}(V) + 2^{n_{j+k-1}(V)}$$
$$= N_{0,j_0}(V) + 2^{n_{j_0+k-1}(V)}$$

Ce qui implique :

$$\frac{N_{0,j_0+k}(V)}{2^{n_{j_0+k}(V)}} = \frac{1}{2}\frac{N_{0,j_0+k-1}(V)}{2^{n_{j_0+k}(V)}} + \frac{1}{2} = \frac{1}{2}\frac{N_{0,j}(V)}{2^{n_{j_0+k}(V)}} + \frac{1}{2}$$

Quelque soit la valeur de réel positif $\varepsilon$ on a :

$$\frac{1}{2} < \frac{N_{0,j_0+k}(V)}{2^{n_{j_0+k}(V)}} < \frac{\varepsilon}{2} + \frac{1}{2}$$

Si on prend le cas ou $\varepsilon < \frac{1}{2}$, on déduit que :

$$\varepsilon < \frac{1}{2} < \frac{N_{0,j_0+k}(V)}{2^{n_{j_0+k}(V)}}$$

Absurde parce que on a :

$$\forall\, \varepsilon > 0, \exists\, j_0 \text{ tel que } \forall\, j > j_0 \text{ on a } \frac{N_{0,j}(V)}{2^{n_j(V)}} < \varepsilon$$

Et comme $j_{0+k} > j_0$ donc nécessairement il faut que :

$$\frac{N_{0,j_0+k}(V)}{2^{n_{j_0+k}(V)}} < \varepsilon$$

On peut conclure que lorsque $N_{0,j}(V)$ tend vers l'infini, le rapport $N_{0,j}(V)/2^{n_j(V)}$ ne peut pas tendre vers pas vers 0. Par contre l'autre cas est évident si $N_{0,j}(V)$ tend vers un entier naturel bien déterminé donc $N_{0,j}(V)/2^{n_j(V)}$ tend vers 0.

**Corollaire 6.1.** La suite $(r_{j,0})_{j\geq 1}$ tend vers 0 lorsque j tend vers l'infini **si et seulement si** la suite $(N_{0,j})_{j\geq 1}$ tend vers une limite fini lorsque j tend vers l'infini. Autrement dit:

$$\lim_{j\to+\infty} \frac{N_{0,j}(V)}{2^{n_j(V)}} = 0 \iff \lim_{j\to+\infty} N_{0,j}(V) \text{ est finie} \qquad (6.4)$$

**Théorème 6.2.** Soit V un vecteur de parité de **longueur infinie** donc V est **réalisable ou convertible en une suite de Collatz si et seulement si** la condition suivante est vérifiée :

$$\lim_{j\to+\infty} \frac{N_{0,j}(V)}{2^{n_j(V)}} = 0$$

Le vecteur binaire infinie V est non réalisable ou non convertible en suite de Collatz **si et seulement si** V satisfait la condition suivante :

$$\lim_{j \to +\infty} \frac{N_{0,j}(V)}{2^{n_j(V)}} \neq 0$$

**Théorème 6.3.** Soit V un vecteur de parité de longueur infinie rappeler que :

$$q_j(V) = \frac{X_j(V)}{2^{n_j(V)}}$$

On considère la suite caractéristique $(q_j(V))_{j \geq 1}$ donc V est réalisable **si et seulement si** la suite $(X_j(V)/2^{n_j(V)})_{j \in \mathbb{N}^*}$ est équivalente à une suite des entiers naturels non nuls c'est à dire que $q_j(V)$ est extrêmement proche d'un entier naturel non nul quelconque pour des valeurs de j très grandes autrement la condition suivante est vérifiée :

$$V \in \mathbb{B}_R \Leftrightarrow \frac{X_j(V)}{2^{n_j(V)}} \sim k_j \in \mathbb{N}^* \tag{6.5}$$

Equivaut à :

$$\frac{P_j(V) a_j(V)}{2^{n_j(V)}} \sim k_j \in \mathbb{N}^* \tag{6.6}$$

Par opposition, V est **non réalisable ou non convertible en une suite de Collatz** si et seulement si la condition suivante est vérifiée :

$$V \notin \mathbb{B}_R \Leftrightarrow \frac{X_j(V)}{2^{n_j(V)}} \nsim k_j \in \mathbb{N}^* \tag{6.7}$$

**Corollaire 6.2.** Soit V un vecteur de parité de longueur infinie alors V est convertible en suite de Collatz si la condition suivante est remplie :

$$V \in \mathbb{B}_R \Leftrightarrow \frac{X_j^*(V)}{2^{n_j(V)}} \sim K_j^* \in \mathbb{N}^* \tag{6.8}$$

**Remarque 6.1.** Selon la conjecture de Collatz, on a :

$$\mathbb{B}_R \subset \mathbb{B}_{(0,1)}$$

Ce qui nous permet de conclure que selon la conjecture de Collatz on a :

$$\frac{X_j(V)}{2^{n_j(V)}} \sim k_j \Rightarrow V \in \mathbb{B}_{(0,1)} \tag{6.9}$$

**Autres conditions pour qu'un vecteur de parité infini soit réalisable ou convertible en une suite infinie suite de Collatz**

On sait que :

$$P_j(V) = 2^{n_j(V)} A_j(V) + B_j(V) \text{ avec } B_j(V) < 2^{n_j(V)}$$

Ce qui nous permet d'écrire :

$$\frac{P_j(V)\, a_j(V)}{2^{n_j(V)}} = A_j(V) a_j(V) + \frac{B_j(V)\, a_j(V)}{2^{n_j(V)}}$$

On suppose que :

$$B_j(V)\, a_j(V) = 2^{n_j(V)}\, f_{j,1}(V) + f_{j,2}(V)$$

Tel que $f_{j,1}(V)$ et $f_{j,2}(V)$ sont deux entiers naturels non nuls.

Donc on peut conclure :

$$\frac{X_j(V)}{2^{n_j(V)}} \sim k_j \Leftrightarrow \left( \lim_{j \to +\infty} \frac{f_{j,2}(V)}{2^{n_j(V)}} = 0 \text{ ou bien } \lim_{j \to +\infty} \frac{f_{j,2}(V)}{2^{n_j(V)}} = 1 \right) \qquad (6.10)$$

On peut conclure qu'une étude approfondie sur les nombres caractéristique et les équations caractéristiques relatives aux vecteurs de parité de longueur infinies est nécessaire pour déterminer pour quelles conditions, les valeurs de $X_j(V) = P_j(V)\, a_j(V)$ est très voisine d'un multiple de $2^{n_j(V)}$. Selon **la conjecture de Collatz** les seules conditions sont :

$$\lim_{j \to +\infty} \frac{m_j(V)}{n_j(V)} = \frac{1}{2}$$

$$\lim_{j \to +\infty} \frac{P_j(V)}{2^{n_j(V)}} = 1$$

En effet, ces deux limites correspondent aux vecteurs appartenant à l'ensemble $\mathbb{B}_{(0,1)}$

### Références